\documentclass[10pt,journal,final,twocolumn,twoside,letterpaper]{IEEEtran}
\usepackage{graphicx,calc}
\usepackage{epstopdf}
\usepackage[latin1]{inputenc}
\usepackage{multirow}
\usepackage{amsmath,amsthm,amssymb,amscd,stmaryrd}
\usepackage{arydshln}
\def\bestapprox{\sigma}
\def\pdf{p}
\def\cdf{F}
\def\slevel{\rho}
\def\mlevel{\delta}
\newcommand\Gfun[2]{G_{#2}[#1]}
\newcommand\Hfun[1]{H[#1]}

\def\Expect{\mathbb{E}}

\def\x{\mathbf{x}}
\def\n{\mathbf{n}}
\def\y{\mathbf{y}}

\def\R{\mathbb{R}}

\DeclareMathOperator*{\sgn}{sgn}

\def\z{\mathbf{z}}
\def\Prob{\mathbb{P}}

\def\sensing{\mathbf{\Phi}}
\def\Ephi{\mathbb{E}_{\sensing}}

\newtheorem{prop}{Proposition}
\newtheorem{remark}{Remark}
\newtheorem{lemma}{Lemma}
\newtheorem{theorem}{Theorem}
\newtheorem{example}{Example}
\newtheorem{corollary}{Corollary}
\newtheorem{definition}{Definition}

\DeclareMathOperator*{\argmin}{argmin}
\begin{document}


\title{Compressible Distributions\\ for High-dimensional Statistics}

\author{R\'emi Gribonval, {\em Senior Member, IEEE},
Volkan Cevher, {\em Senior Member, IEEE},\protect\\ and
Mike E. Davies, {\em Senior Member, IEEE}
\thanks{R{\'e}mi Gribonval is with INRIA, Centre Inria Rennes - Bretagne Atlantique, 35042 Rennes Cedex, France.}
\thanks{Volkan Cevher is with the Laboratory for Information and Inference
Systems (LIONS),  Ecole Polytechnique Federale
de Lausanne, Switzerland.}
\thanks{Mike Davies is with the Institute for Digital Communications
  (IDCOM) \& Joint Research Institute for Signal and Image Processing,
  School of Engineering and Electronics, University of Edinburgh, The
  King's Buildings, Mayfield Road, 
Edinburgh EH9 3JL, UK.}
\thanks{This work was supported in part by the European Commission through the project SMALL (Sparse Models, Algorithms and Learning for Large-Scale data) under FET-Open, grant number 225913.}
\thanks{This work was supported in part by part by the European Commission under Grant MIRG-268398, ERC Future Proof 279735,  and DARPA KeCoM program \#11-DARPA-1055. VC also would like to acknowledge Rice University for his Faculty Fellowship.}
\thanks{MED acknowledges support of his position from the Scottish Funding Council and their support of the Joint Research Institute in Signal and Image Processing with the Heriot-Watt University as a component of the Edinburgh Research Partnership. This work was supported in part by the UK Engineering and Physical Science Research Council, grant EP/F039697/1.}
\thanks{Copyright (c) 2012 IEEE. Personal use of this material is permitted.  However, permission to use this material for any other purposes must be obtained from the IEEE by sending a request to pubs-permissions@ieee.org.}
}

\markboth{to appear in IEEE Transactions on Information Theory, 2012}{R. Gribonval, V. Cevher and M.E. Davies: Compressible Distributions}
%
%
%
%
%
%

\maketitle
\IEEEpeerreviewmaketitle

\begin{abstract}
We develop a principled way of identifying probability distributions
whose independent and identically distributed realizations are
compressible, i.e., can be well-approximated as sparse. We focus on
Gaussian compressed sensing, an example of underdetermined linear regression, 
where compressibility is known to ensure the success of estimators
exploiting sparse regularization. We prove that many distributions
revolving around maximum a posteriori (MAP) interpretation of sparse
regularized estimators are in fact incompressible, in the limit of
large problem sizes. We especially highlight the Laplace distribution and
$\ell^{1}$ regularized estimators such as the Lasso and Basis Pursuit
denoising. We rigorously disprove the myth that the success of $\ell^1$
minimization for compressed sensing image reconstruction is a simple
corollary of a Laplace model of images combined with Bayesian MAP
estimation, and show that in fact quite the reverse is true.
To establish this result, we identify non-trivial undersampling
regions where the simple least squares solution almost surely
outperforms an oracle sparse solution, when the data is generated
from the Laplace distribution. We also provide simple rules of thumb to
characterize classes of compressible and incompressible
distributions based on their second and fourth moments. Generalized
Gaussian and generalized Pareto distributions serve as running
examples. 
\end{abstract} 




\IEEEpeerreviewmaketitle

\section{Introduction}\label{cp: intro}
High-dimensional data is shaping the current {\em modus operandi} of statistics. Surprisingly, while the ambient dimension is large in many problems, natural constraints and parameterizations often cause data to cluster along low-dimensional structures. Identifying and exploiting such structures using probabilistic models is therefore quite important for statistical analysis, inference, and decision making.

In this paper, we discuss compressible distributions, whose independent and identically distributed (iid) realizations can be well-approximated as sparse.
%
Whether or not a distribution is compressible is important in the
context of many applications, among which we highlight two here:
statistics of natural images, and statistical regression for linear
inverse problems such as those arising in the context of compressed sensing.
\subsubsection*{Statistics of natural images} 
Acquisition, compression, denoising, and analysis of natural images (similarly, medical, seismic, and hyperspectral images) draw high scientific and commercial interest. 
Research  to date in natural image modeling has had two distinct approaches, with one focusing on deterministic explanations and the other pursuing probabilistic models. Deterministic approaches (see e.g.~\cite{Chang2000,Choi:1999aa}) operate under the assumption that the transform domain representations (e.g., wavelets, Fourier, curvelets, etc.) of images are ``compressible''. Therefore, these approaches threshold the transform domain coefficients for sparse approximation, which can be used for compression or denoising.

Existing probabilistic approaches also exploit coefficient decay in transform domain representations, and learn probabilistic models by approximating the coefficient histograms or moment matching. For natural images, the canonical approach (see e.g.~\cite{portilla2003image}) is to fit probability density functions (PDF's), such as generalized Gaussian distributions and the Gaussian scale mixtures, to the histograms of wavelet coefficients while trying to simultaneously capture the dependencies observed in their marginal and joint distributions.

\subsubsection*{Statistical regression} Underdetermined linear
regression is a fundamental problem in statistics, applied
mathematics, and theoretical computer science with broad applications---from subset selection to compressive
sensing~\cite{MR2241189,Candes:2004aa}  and inverse problems (e.g.,
deblurring), and from data streaming to error corrective coding. In
each case, we seek an unknown vector $\x \in \R^N$, given its dimensionality
reducing, linear projection $\y \in \R^m$ ($m < N$) obtained via a
known encoding matrix $\sensing \in \R^{m\times N}$, as 
\begin{equation}
 \y = \sensing \x + \n,
\end{equation}
where $\n \in \R^m$ accounts for the perturbations in the linear system, such as physical noise.
The core challenge in decoding $\x$ from $\y$ stems from the simple fact that dimensionality reduction loses information in general: for any vector $v\in \text{kernel}(\sensing)$, it is impossible to distinguish $\x$ from $\x+v$ based on $\y$ alone. 

Prior information on $\x$ is therefore necessary to estimate the true $\x$ among the infinitely many possible solutions.
It is now well-known that geometric sparsity models (associated to approximation of $\x$ from a finite union of low-dimensional subspaces in $\R^N$ \cite{DBLP:journals/tit/BlumensathD09}) play an important role in obtaining ``good'' solutions. A widely exploited decoder is the $\ell^{1}$ decoder $\Delta_{1}(\y) := \arg\min_{\tilde{\x} : \y = \sensing \tilde{\x}} \|\tilde{\x}\|_{1}$ whose performance can be explained via the geometry of projections of the $\ell^{1}$ ball in high dimensions~\cite{Donoho:2009aa}.
A more probabilistic perspective considers $\x$ as drawn from a 
  distribution. As we will see, compressible iid
distributions~\cite{Baraniuk:2010aa,IEEE_NIPS09_Cevher} countervail
the ill-posed nature of compressed sensing problems by generating
vectors that, in high dimensions, are well approximated by the geometric sparsity model.

\subsection{Sparsity, compressibility and compressible distributions}

A celebrated result from compressed sensing
\cite{MR2241189,stablesigrecovery-CandesRombergTao-2006} 
is that under certain conditions, a $k$-{\em sparse} vector $\x$ (with
only $k$ non-zero entries where $k$ is usually much smaller than $N$) can be exactly recovered from its noiseless
projection $\y$ using the $\ell^1$
decoder, as long as $m \gtrsim k \log N/k$. Possibly the most
striking result of this type is the Donoho-Tanner weak phase transition
that, for Gaussian sensing matrices, completely characterizes the typical success or failure of the
$\ell^1$ decoder in the large scale limit \cite{Donoho:2009aa}. 


Even when the vector $\x$ is not sparse, under certain ``compressibility'' conditions typically expressed in terms of (weak) $\ell^{p}$ balls, the
$\ell^{1}$-decoder provides estimates with controlled
accuracy~\cite{Cohen:2006aa,Candes:2008aa,Davies:2009aa,Donoho:2011}. 
Intuitively one should only expect a sparsity-seeking estimator to
perform well if the vector being reconstructed is at least
approximately sparse. Informally,
compressible vectors can be defined as follows: 
\begin{definition}[Compressible vectors]\label{def: relapprox} 
 Define the {\em relative} best $k$-term approximation error $\bar{\bestapprox}_k(\x)_q$ of a vector $\x$ as
\begin{equation}\label{eq: rel app err}
  \bar{\bestapprox}_k(\x)_q= \frac{{\bestapprox}_k(\x)_q}{\|\x\|_q},
\end{equation}
where $\bestapprox_k(\x)_q := \inf_{\|\y\|_0 \leq k} \|\x-\y\|_q$  is the best $k$-term approximation error of $\x$, and $\|\x\|_q$ is the $\ell^q$-norm of $\x$, $q\in(0,\infty)$. By convention $\|{\x}\|_0$ counts the non-zero coefficients of ${\x}$.
A vector $\x \in \R^{N}$ is $q$-{\em compressible} if $ \bar{\bestapprox}_k(\x)_q \ll 1$ for some $k \ll N$.
\end{definition}

This definition of compressibility differs slightly from those that
are closely linked to weak $\ell^p$ balls in that, above, we consider
{\em relative} error. This is discussed further in Section~\ref{sec:
  instance optimality}. 

When moving from the deterministic setting to the stochastic setting
it is natural to ask when reconstruction guarantees equivalent to the deterministic
ones exist. The case of typically sparse vectors is most easily dealt
with and can be characterized by a distribution with a probability
mass of $(1-k/N)$ at zero, e.g., a Bernoulli-Gaussian distribution. Here the results of Donoho and Tanner still
apply as a random vector drawn from such a distribution is typically sparse, with approximately $k$ nonzero entries, while the $\ell^1$ decoder is blind to
the specific non-zero values of $\x$.

The case of compressible vectors is less straightforward: when is a vector
generated from iid draws of a given 
distribution typically compressible? This is the question
investigated in this paper. To exclude the sparse case, we restrict
ourselves to distributions with a well defined density $\pdf(x)$.

Broadly speaking, we can define compressible distributions as follows. 
\begin{definition}[Compressible distributions]\label{def: cp}
 Let $X_n (n \in \mathbb{N})$ be iid samples from a probability distribution with probability density function (PDF) $\pdf(x)$, and $\x_N = (X_1,\ldots,X_N) \in \R^N$.
The PDF $\pdf(x)$ is said to be $q$-{\em compressible} with parameters $(\epsilon,\kappa)$ when
  \begin{equation}
\limsup_{N \to \infty} \bar{\bestapprox}_{k_N}(\x_N)_q
\stackrel{a.s.}{\le} \epsilon, (\text{a.s.: almost surely});
\end{equation}
for any sequence $k_N$ such that $\liminf_{N \to \infty} \frac{k_N}{N} \ge \kappa$.

\end{definition}
The case of interest is when $\epsilon\ll 1$ and $\kappa \ll 1$: iid
realizations of a $q$-compressible distribution with parameters
$(\epsilon,\kappa)$ live in $\epsilon$-proximity to the union of
$\kappa N$-dimensional hyperplanes, where the closeness is measured in
the $\ell^q$-norm. These hyperplanes are aligned with the coordinate
axes in $N$-dimensions. 

One can similarly define an incompressible distribution as:
\begin{definition}[Incompressible distributions]\label{def: incp}
 Let $X_n$ and $\x_N$ be defined as above. The PDF
 $\pdf(x)$ is said to be $q$-{\em incompressible} with parameters
 $(\epsilon,\kappa)$ when 
  \begin{equation}
\liminf_{N \to \infty} \bar{\bestapprox}_{k_N}(\x_N)_q
\stackrel{a.s.}{\ge} \epsilon, 
\end{equation}
for any sequence $k_N$ such that $\limsup_{N \to \infty} \frac{k_N}{N} \le \kappa$.
\end{definition}
This states that the iid realizations of an incompressible
distribution live away from the $\epsilon$-proximity of the union of
$\kappa N$-dimensional hyperplanes, where $\epsilon \approx 1$. 

More formal characterizations of the  ``compressibility'' or the
``incompressibility'' of a distribution with PDF $\pdf(x)$ are investigated in
this paper. With a special emphasis on the context of compressed
sensing with a Gaussian encoder $\sensing$, we discuss and characterize the compatibility of such distributions with extreme levels of undersampling. 
As a result, our work features both positive and negative conclusions on achievable approximation performance of probabilistic modeling in compressed sensing\footnote{Similar ideas were recently proposed in \cite{Amini:2011}, however, while
the authors explore the stochastic concepts of compressibility they do not examine the
implications for signal reconstruction in compressed sensing type
scenarios.}.

\subsection{Structure of the paper}
The main results are stated in Section~\ref{sec:detailedoverview} together with a discussion of their conceptual implications. The section is concluded by Table~\ref{table: summary}, which provides an overview at a glance of the results. The following sections discuss in more details our contributions, while the bulk of the technical contributions is gathered in an appendix, to allow the main body of the paper to concentrate on the conceptual implications of the results. 
As running examples, we focus on the Laplace distribution for incompressibility and the generalized Pareto distribution for compressibility, with a Gaussian encoder $\sensing$. 


\section{Main results}
\label{sec:detailedoverview}
%
%
%
%
%
%

In this paper, we aim at bringing together the deterministic and probabilistic models of compressibility in a simple and general manner under the umbrella of compressible distributions. To achieve our goal, we dovetail the concept of order statistics from probability theory with the deterministic models of compressibility from approximation theory. 

Our five ``take home'' messages for compressed sensing are as follows:
\begin{enumerate}
\item $\ell^1$ minimization does not assume that the underlying
  coefficients have a Laplace distribution. In fact, the relatively flat
  nature of vectors drawn iid from a Laplace distribution makes them, in some sense, the worst for
  compressed sensing problems. 
\item It is simply not true that the success of $\ell^1$ minimization
  for compressed sensing reconstruction is a simple corollary of
  a Laplace model of data coefficients combined with Bayesian MAP estimation, in
  fact quite the reverse. 
\item Even with the strongest possible recovery guarantees~\cite{Cohen:2006aa,Davies:2009aa}, compressed sensing reconstruction of Laplace
  distributed vectors with the $\ell^{1}$ decoder offers no guarantees beyond the trivial
  estimator, $\widehat{\x} = 0$. 
\item More generally, for high-dimensional vectors $\x$ drawn iid from any density with
  bounded fourth moment $\Expect X^{4} < \infty$, even with the help of a sparse oracle, there
  is a critical level of undersampling below which the sparse oracle
  estimator is worse (in relative $\ell^2$ error) than the simple least-squares
  estimator.
\item In contrast, when a high-dimensional vector $\x$ is drawn from a density with
  infinite second moment $\Expect X^{2} = \infty$, then the $\ell^1$ decoder can reconstruct
  $\x$ with arbitrarily small relative $\ell^2$ error.
\end{enumerate}

\subsection{Relative sparse approximation error} 
By using Wald's lemma on order statistics, we characterize the relative sparse approximation errors of iid PDF realizations, whereby providing solid mathematical ground to the earlier work of Cevher~\cite{IEEE_NIPS09_Cevher} on compressible distributions.
While Cevher exploits the decay of the expected order statistics, his approach is inconclusive in characterizing the ``incompressibility'' of distributions. We close this gap 
by introducing a function $\Gfun{\pdf}{q}(\kappa)$ so that iid vectors as in Definition~\ref{def: cp} satisfy
\(
\lim_{N \to \infty} \bar{\bestapprox}_{k_N}(\x_N)_q^q
\stackrel{a.s.}{=} \Gfun{\pdf}{q}(\kappa)
\)
when $\lim_{N \to \infty} k_N/N = \kappa \in (0,1)$.

\begin{prop}
\label{prop:RelativeError}
Suppose $\x_N \in \R^N$ is iid with respect to $\pdf(x)$ as in Definition \ref{def: cp}. Denote $\bar{\pdf}(x) := 0$ for $x<0$, and $\bar{\pdf}(x) := \pdf(x) + \pdf(-x)$ for $x\geq 0$ as the PDF of $|X_n|$, and $\bar{\cdf}(t) := \mathbb{P}(|X|\leq t)$ as its cumulative density function (CDF). Assume that $\bar{\cdf}$ is continuous and strictly increasing on some interval $[a\ b]$, with $\bar{\cdf}(a) = 0$ and $\bar{\cdf}(b) = 1$, where $0 \leq a < b \leq \infty$.  For any $0<\kappa\leq1$, define the  following function:
\begin{equation}
\label{eq:DefGfun}
\Gfun{\pdf}{q}(\kappa) := \frac{\int_0^{\bar{\cdf}^{-1}(1-\kappa)} x^q
 \bar{\pdf}(x) dx}{\int_0^\infty x^q \bar{\pdf}(x) dx}.
\end{equation}

\begin{enumerate}
\item {\textbf{Bounded moments:}} assume $\Expect |X|^q < \infty$ for some $q\in (0,\infty)$. Then, $\Gfun{\pdf}{q}(\kappa)$ is also well defined for $\kappa=0$, and for any sequence $k_N$ such that $\lim_{N \to \infty} \frac{k_N}{N} = \kappa \in [0,1]$, the following holds almost surely
\begin{equation}
\lim_{N \to \infty} \bar{\bestapprox}_{k_N}(\x_N)_q^q
\stackrel{a.s.}{=} \Gfun{\pdf}{q}(\kappa).
\end{equation}
\item {\textbf{Unbounded moments:}} assume $\Expect |X|^q = \infty$ for some $q\in (0,\infty)$. Then, for $0 < \kappa \leq 1$ and any sequence $k_N$ such that $\lim_{N \to \infty} \frac{k_N}{N} = \kappa$, the following holds almost surely
\begin{equation}
\lim_{N \to \infty} \bar{\bestapprox}_{k_N}(\x_N)_q^q
\stackrel{a.s.}{=} \Gfun{\pdf}{q}(\kappa) = 0.
\end{equation}
\end{enumerate}
\end{prop}

Proposition \ref{prop:RelativeError} provides a principled way of obtaining the compressibility parameters $(\epsilon,\kappa)$ of distributions in the high dimensional scaling of the vectors. An immediate application is the incompressibility of the Laplace distribution.

\begin{example}
\label{ex:Laplace}
As a stylized example, consider the Laplace distribution (also known as the double exponential) with scale parameter 1, whose PDF is given by
\begin{equation}
\label{eq:DefLaplace}
\pdf_1(x) := \frac{1}{2} \exp(-|x|).
\end{equation}
We compute in Appendix~\ref{app:Laplace}:
\begin{align}
\label{eq:BestKTermLaplaceAsympL1}
\Gfun{\pdf_1}{1}(\kappa) &= 1-\kappa \cdot \Big(1 + \ln 1/\kappa \Big),\\
\label{eq:BestKTermLaplaceAsymp}
\Gfun{\pdf_1}{2}(\kappa)
&=1-\kappa \cdot \Big(1 + \ln 1/\kappa + \frac12 (\ln 1/\kappa)^2\Big).
\end{align}
Therefore, it is straightforward to see that the Laplace distribution is {\em not $q$-compressible} for $q\in \{1,2\}$: it is not possible to simultaneously have both $\kappa$ and $\epsilon = \Gfun{\pdf_1}{q}(\kappa)$ small.
\end{example} 

\subsection{Sparse modeling vs.\ sparsity promotion} We show that the {\em maximum a posteriori} (MAP) interpretation of standard deterministic sparse recovery algorithms is, in some sense, inconsistent.  To explain why, we consider the following decoding approaches to estimate a vector $\x$ from its encoding $\y = \sensing \x$:
\begin{align}\label{eq:P1}
\Delta_1(\y) & =\argmin_{\tilde{\x}:\y = \sensing \tilde{\x}} \|\tilde{\x}\|_1,\\
\label{eq:DefLSSolution}
\Delta_{\textrm{LS}}(\y)  &= \argmin_{\tilde{\x}:\y = \sensing \tilde{\x}}\|\tilde{\x}\|_2 = \sensing^+ \y,\\
\label{eq:DefOSSolution}\Delta_{\textrm{oracle}}(\y,\Lambda) &= \argmin_{\tilde{\x}:\textrm{support}(\tilde{\x}) = \Lambda}\|\y-\sensing\tilde{\x}\|_2 = \sensing_{\Lambda}^+ \y,\\
\label{eq:DeftrSolution}\Delta_{\textrm{trivial}}(\y) &= 0.
\end{align}
Here, $\sensing_\Lambda$ denotes the sub-matrix of $\sensing$
restricted to the columns indexed by the set $\Lambda$. The decoder $\Delta_1$ regularizes the solution space via the $\ell^1$-norm. It is the {\em de facto} standard Basis Pursuit formulation \cite{CD99} for sparse recovery, and is tightly related to the Basis Pursuit denoising  (BPDN) and the least absolute shrinkage and selection operator (LASSO) \cite{Tibshirani94regressionshrinkage}:
\begin{align*}
 \Delta_{\text{BPDN}}(\y)&=  \argmin_{\tilde{\x}}  \left\{\frac{1}{2} \|\y - \sensing \tilde{\x}\|_2^2 + \lambda \|\tilde{\x}\|_1\right\}
\end{align*}
where $\lambda$ is a constant. Both $\Delta_1$ and the BPDN formulations can be solved in polynomial time through convex optimization techniques.
The decoder $\Delta_{\textrm{LS}}$ is the traditional minimum least-squares solution, which is related to the Tikhonov regularization or ridge regression.
It uses the Moore-Penrose pseudo-inverse $\sensing^+ = \sensing^T(\sensing\sensing^T)^{-1}$. 
The oracle sparse decoder $\Delta_{\textrm{oracle}}$ can be seen as an idealization of sparse decoders, which combine subset selection (the choice of $\Lambda$) with a form of linear regression. It is an ``informed'' decoder that has the side information of the
index set $\Lambda$ associated with the largest components in $\x$.
The trivial decoder $\Delta_{\textrm{trivial}}$  plays the devil's advocate for the performance guarantees of the other  decoders.

\subsubsection{Almost sure performance of decoders}

When the encoder $\sensing$ provides near isometry to the set of sparse vectors \cite{stablesigrecovery-CandesRombergTao-2006}, the decoder $\Delta_1$ features an {\em instance optimality} property \cite{Cohen:2006aa,Davies:2009aa}:
\begin{equation}
\label{eq:DefInstanceOptimality}
\|\Delta_1(\sensing \x)-\x\|_1 \leq C_k(\sensing) \cdot \bestapprox_k(\x)_1, \forall \x;
\end{equation}
where $C_k(\sensing)$ is a constant which depends on $\sensing$. 
A similar result holds with the $\|\cdot\|_2$ norm on the left hand side. Unfortunately, it is impossible to have the same uniform guarantee for all $\x$ with $\bestapprox_k(\x)_2$ on the right hand side \cite{Cohen:2006aa}, but for any given $\x$, it becomes possible {\em in probability} \cite{Cohen:2006aa,DeVore2009}. 
For a Gaussian encoder, $\Delta_1$ recovers exact sparse vectors perfectly from as few as
$m \approx 2e k \log N/k$ with high probability \cite{Donoho:2009aa}.
\begin{definition}[Gaussian encoder]\label{def: GaussianEncoder}
 Let $\phi_{i,j}$, $i,j \in \mathbb{N}$ be iid Gaussian variables
 $\mathcal{N}(0,1)$. The $m \times N$ Gaussian encoder is the random matrix $\sensing_N := \left[\phi_{ij}/\sqrt{m} \right]_{1 \leq i \leq m, 1\leq j \leq N}$.
\end{definition}
{\em In the sequel, we only consider the Gaussian encoder},
leading to Gaussian compressed sensing (G-CS) problems. 
In Section~\ref{sec:expectedperf}, we theoretically characterize the almost sure performance of the estimators $\Delta_{\textrm{LS}}$, $\Delta_{\textrm{oracle}}$ for arbitrary high-dimensional vectors $\x$. 
We concentrate our analysis to the noiseless setting\footnote{Coping with noise in such problems is important both from a practical and a statistical perspective. Yet, the noiseless setting is relevant to establish negative results such as Theorem~\ref{Th:4MomentIntro} which shows the {\em failure} of sparse estimators {\em in the absence of noise}, for an 'undersampling ratio' $\delta$ bounded away from zero. Straightforward extensions of more positive results such as Theorem~\ref{th:InfiniteMomentInstOpt} to the Gaussian noise setting can be envisioned.} ($\mathbf{n} = 0$). 
The least squares decoder $\Delta_{\textrm{LS}}$ has expected performance $\Ephi \|\Delta_{\textrm{LS}}(\sensing\x)-\x\|_{2}^{2}/\|\x\|_{2}^{2} = 1-\mlevel$, independent of the vector $\x$, where 
\begin{equation}
\label{eq:DefUndersamplingRatio}
\mlevel := m/N
\end{equation}
is the {\em undersampling ratio} associated to the matrix $\sensing$ (this terminology comes from compressive sensing, where $\sensing$ is a sampling matrix).
In theorem~\ref{th:ExpectedPerfOracle} the expected performance of the
oracle sparse decoder $\Delta_{\textrm{oracle}}$  is shown to satisfy
\[\frac{\Ephi \|\Delta_{\textrm{oracle}}(\sensing\x,\Lambda)-\x\|_2^2}{ \|\x\|_2^2} = \frac{1}{1-\frac{k}{m-1}} \times
\frac{\bestapprox_k(\x)^2}{\|\x\|_2^2}.\]
This error is the balance between two factors. The first factor grows
with $k$ (the size of the set $\Lambda$ of largest entries of $\x$
used in the decoder) and  reflects the (ill-)conditioning of the
Gaussian submatrix $\sensing_\Lambda$.
The second factor is the best $k$-term relative approximation error, which shrinks as $k$ increases.
This highlights the inherent trade-off present in any sparse estimator, namely the level of sparsity $k$
versus the conditioning of the sub-matrices of $\sensing$.

\subsubsection{A few surprises regarding sparse recovery guarantees}
We highlight two counter-intuitive results below:

\paragraph{A crucial weakness in appealing to instance optimality}

Although instance optimality~\eqref{eq:DefInstanceOptimality} is usually considered as a strong  property, it involves an implicit trade off: when $k$ is small, the $k$-term error $\bestapprox_k(\x)$ is large, while for larger $k$, the constant $C_k(\sensing)$ is large. For instance, we have $C_k(\sensing) = \infty$, when $k \geq m$.  

In Section~\ref{sec: instance optimality} we provide new key insights for instance optimality of algorithms.
Informally, we show that when $\x_{N} \in \R^N$ is iid with respect to $\pdf(x)$ as in Definition \ref{def: cp}, and when $\pdf(x)$ satisfies the hypotheses of Proposition~\ref{prop:RelativeError}, if
\begin{equation}
\label{eq:IOfails}
 \Gfun{\pdf}{1}(\kappa_0) \geq 1/2,
\end{equation}
where $\kappa_0\approx 0.18$ is an absolute constant, then the best possible upper bound in the instance optimality~\eqref{eq:DefInstanceOptimality} for a Gaussian encoder satisfies (in the limit of large $N$)
\[
 C_{k}(\sensing) \cdot \bestapprox_{k}(\x)_{1} \geq  \|\x\|_{1} = \|\Delta_{\textrm{trivial}}(\x)-\x\|_{1}.
\]
In other words, {\em for distributions with PDF $\pdf(x)$ satisfying~\eqref{eq:IOfails}, in high dimension $N$, 
instance optimality results for the decoder $\Delta_1$ with a Gaussian encoder can at best guarantee the performance (in the $\ell^{1}$ norm) of the trivial decoder $\Delta_{\textrm{trivial}}$!}

Condition \eqref{eq:IOfails} holds true for many general PDF's;  it is easily verifiable for the Laplace distribution based on Example~\ref{ex:Laplace}, and explains the observed failure of the $\ell^{1}$ decoder on Laplace data~\cite{Seeger2008}. This is  discussed further in Section~\ref{sec: instance optimality}.

\paragraph{Fundamental limits of sparsity promoting decoders}
The expected $\ell^{2}$ relative error of the least-squares estimator $\Delta_{\textrm{LS}}$ degrades linearly as $1-\mlevel$ with the undersampling factor $\mlevel := m/N$, and
therefore does not provide good reconstruction at low sampling rates $\mlevel \ll 1$. It is therefore quite surprising that we can determine a large class of distributions for which the oracle sparse decoder $\Delta_{\textrm{oracle}}$ is outperformed by the simple least-squares decoder $\Delta_{\textrm{LS}}$.
\begin{theorem}
\label{Th:4MomentIntro}
Suppose that $\x_{N} \in \R^N$ is iid with respect to $\pdf(x)$ as in Definition \ref{def: cp}, and that $\pdf(x)$ satisfies the hypotheses of Proposition~\ref{prop:RelativeError} and has a finite fourth-moment
\[
\Expect X^4 < \infty.\]
There exists a minimum undersampling ratio $\mlevel_0$ with the following property: for any $\slevel \in (0,1)$, if $\sensing_{N}$ is a sequence of $m_{N} \times N$ Gaussian encoders with $\lim_{N \to \infty} m_{N}/N = \mlevel <\mlevel_0$, and $\lim_{N \to \infty} k_{N}/m_{N} = \slevel$, then we have almost surely
\begin{align*}
\lim_{N \to \infty} & 
\frac{\|\Delta_{\textrm{oracle}}(\sensing_N \x_{N},\Lambda_N)-\x_N\|_2^2}{ \|\x_N\|_2^2}\\
&\stackrel{a.s.}{=} \frac{\Gfun{\pdf}{2}(\slevel\mlevel)}{1-\slevel}\\
& > 1-\delta \stackrel{a.s.}{=}
\lim_{N \to \infty}
\frac{\|\Delta_{\textrm{LS}}(\sensing_N \x_{N})-\x_N\|_2^2}{ \|\x_N\|_2^2}.
\end{align*}
\end{theorem}
Thus if the data PDF $\pdf(x)$ has a finite fourth moment and a continuous CDF, there exists a level of undersampling below which a simple least-squares reconstruction (typically a dense vector estimate) provides an estimate, which is closer to the true vector $\x$ (in the $\ell^{2}$ sense) than oracle sparse estimation!

Section~\ref{sec:ComparisonL2L0} describes how to determine this undersampling boundary, e.g., for the generalized Gaussian distribution. For the Laplace distribution, $\mlevel_0 \approx 0.15$. In other words, when randomly sampling a high-dimensional Laplace vector, it is better to use least-squares reconstruction than minimum $\ell^1$ norm reconstruction (or any other type of sparse estimator), unless the number of measures $m$ is at least $15 \%$ of the original vector dimension $N$. To see how well Theorem \ref{Th:4MomentIntro} is grounded in practice, we provide the following example:

\begin{example}
\label{ex:Laplace2}
Figure~\ref{fig:LaplaceRelError} examines in more detail the
performance of the estimators for Laplace distributed data at various
undersampling values. The horizontal lines indicate  various
signal-to-distortion-ratios (SDR) of $3$dB, $10$dB and $20$dB. Thus
for the oracle estimator to achieve $10$dB, the undersampling rate
must be greater than 0.7, while to achieve a performance level of
$20$dB, something that might reasonably  be expected in many
compressed sensing applications, we can hardly afford any subsampling at all since this requires $\mlevel > 0.9$.

\begin{figure}[htbp]
\begin{center}
\includegraphics[width=9cm]{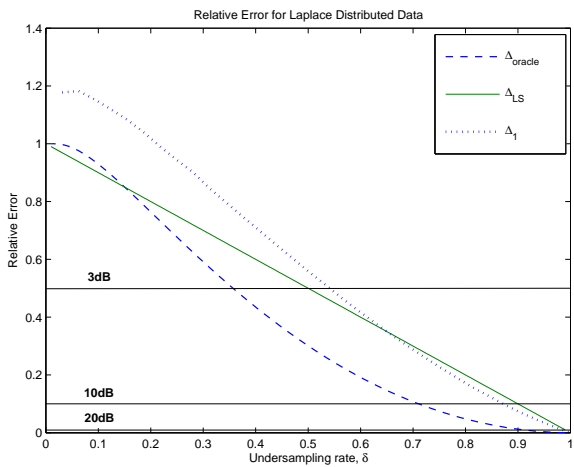}
\caption{\label{fig:LaplaceRelError}
The expected relative error as a function of the undersampling rates $\mlevel$ for data iid from a Laplace distribution using: (a) a linear least squares estimator (solid) and (b) the best oracle sparse estimator (dashed). Also plotted is the empirically observed average relative error  over 5000 instances for the $\Delta_1$ estimator (dotted). The horizontal lines indicate SDR values of $3$dB, $10$dB and $20$dB, as marked.}
\end{center}
\end{figure}
\end{example}

This may come as a shock since, in Bayesian terminology,
$\ell^1$-norm minimization is often conventionally interpreted as the MAP estimator
under the Laplace prior, while least squares is the MAP under the
Gaussian prior. Such MAP interpretations of compressed sensing decoders are further discussed below and contrasted to more geometric interpretations.

\subsection{Pitfalls of MAP ``interpretations'' of decoders}
Bayesian compressed sensing methods employ probability measures as
``priors'' in the 
space of the unknown vector $\x$, and arbitrate the solution space by
using the chosen measure. The decoder $\Delta_{1}$ has a distinct
probabilistic interpretation in the statistics literature. If we presume an iid probabilistic model for $\x$ as $p(X_n)\propto\exp\left(-c|X_n|\right)$ ($n=1,\ldots,N$), then $\Delta_{\text{BPDN}}$ can be viewed as the MAP estimator
\[
\Delta_{\textrm{MAP}}(\y) := \arg\max_{\x} p(\x|\y) = \arg\min_{\x} \{-\log p(\x|\y)\},
\]
when the noise $\n$ is iid Gaussian, which becomes the $\Delta_{1}$
decoder in the zero noise limit. However, as illustrated by Example~\ref{ex:Laplace2}, the decoder $\Delta_{\textrm{MAP}}$ performs quite poorly for iid Laplace vectors. The possible inconsistency of MAP estimators is a known phenomenon~\cite{Nikolova:2007aa}. Yet, the fact that $\Delta_{\textrm{MAP}}$ is outperformed by $\Delta_{\textrm{LS}}$---which is the MAP under the Gaussian prior---when $\x$ is drawn iid according to the Laplacian distribution should remain somewhat counterintuitive to many readers.

It is now not uncommon to stumble upon new proposals in the literature for the modification of $\Delta_1$ or BPDN with diverse thresholding or re-weighting rules based on different hierarchical probabilistic models---many of which  
correspond to a special Bayesian ``sparsity prior'' $\pdf(\x) \propto \exp(-\phi(\x))$ \cite{Wipf08anew}, associated to the minimization of new cost functions
\[
\Delta_{\phi}(\y) := \arg\min_{\x} \frac{1}{2} \|\y-\sensing \x\|_2^2 + \phi(\x).
\]
It has been shown in the context of additive white Gaussian noise denoising that the MAP interpretation of such penalized least-squares regression can be misleading~\cite{GRIBONVAL:2010:INRIA-00486840:1}. Just as illustrated above with $\phi(\x) = \lambda \|\x\|_1$,
 {\em while the geometric interpretations of the cost functions associated to such ``priors'' are useful for sparse recovery, the ``priors'' $\exp(-\phi(\x))$ themselves do not necessarily constitute a relevant ``generative model'' for the vectors.} Hence, such proposals are losing a key strength of the Bayesian approach: the ability to evaluate the ``goodness'' or ``confidence'' of the estimates due to the probabilistic model itself or its conjugate prior mechanics.

In fact, the empirical success of $\Delta_1$ (or $\Delta_{\text{BPDN}}$) results from a combination of two properties: 
\begin{enumerate}
\item the {\em sparsity-inducing} nature of the cost function, due to the non-differentiabi\-lity at zero of the $\ell^1$ cost function; 
\item the {\em compressible} nature of the vector $\x$ to be estimated. 
\end{enumerate}

Geometrically speaking, the objective $\|{\x}\|_{1}$ is related to the $\ell^1$-ball, which intersects with the constraints (e.g., a randomly oriented hyperplane, as defined by $\y=\sensing\x$) along or near the $k$-dimensional hyperplanes ($k\ll N$) that are aligned with the canonical coordinate axes in $\R^N$. The geometric interplay of the objective and the constraints in high-dimensions inherently {\em promotes sparsity}.
An important practical consequence is the ability to design efficient optimization algorithms for large-scale problems, using thresholding operations. Therefore, the decoding process of $\Delta_1$ automatically sifts smaller subsets that best explain the observations, unlike the traditional least-squares $\Delta_{\textrm{LS}}$. 

When $\x_{N}$ has iid coordinates as in Definition~\ref{def: cp}, compressibility is not so much related to the behavior (differentiable or not) of $\pdf(x)$ around zero but rather to the thickness of its tails, e.g., through the necessary property $\Expect X^4 = \infty$ ({\em cf} Theorem~\ref{Th:4MomentIntro}). We further show 
that distributions with infinite variance ($\Expect X^2 = \infty$) almost surely generate vectors which are sufficiently compressible to guarantee that the decoder $\Delta_1$ with a Gaussian encoder $\sensing$ of arbitrary (fixed) small sampling ratio $\mlevel = m/N$ has ideal performance in dimensions $N$ growing to infinity:

\begin{theorem}[Asymptotic performance of the $\ell^1$ decoder under infinite second moment]\label{th:InfiniteMomentInstOpt}
Suppose that $\x_{N} \in \R^{N}$ is iid with respect to $\pdf(x)$ as in Definition~\ref{def: cp}, and that $\pdf(x)$ satisfies the hypotheses of Proposition~\ref{prop:RelativeError} and has infinite second moment
 $\Expect X^2 =
 \infty$.  
Consider a sequence of integers $m_N$ such that $\lim_{N \to \infty}
m_N/N = \mlevel$ where $0<\mlevel<1$ is arbitrary, and let $\sensing_{N}$ be a sequence of $m_{N} \times N$ Gaussian encoders. Then
\begin{equation}\label{eq: Inst OPt l1 performance}
\lim_{N \to \infty}\frac{\| \Delta_1(\sensing_N \x_N)-\x_N \|_2}{\|\x_N\|_2}
\stackrel{a.s.}{=} 0.
\end{equation}
\end{theorem}
%
As shown in Section~\ref{sec:Discussion} there exist PDFs $\pdf(x)$, which combine heavy tails with a non-smooth behavior at zero, such that the associated MAP estimator is sparsity promoting. It is likely that the MAP with such priors can be shown to perform ideally well in the asymptotic regime.



\begin{table*}[ht]\centering
\caption{\label{table: summary}Summary of the main results}
\begin{tabular}{|c||c|c|c|}
  \hline
  Moment property & $\Expect X^2 = \infty$ & $\Expect X^2 < \infty$ and $\Expect X^4 = \infty$  & $\Expect X^4 < \infty$\\
  \hline\hline
  & {\em Theorem~\ref{th:InfiniteMomentInstOpt}}
  & N/A
  & {\em Theorem~\ref{Th:4MomentIntro}}\\
  General result
 &
  $\Delta_1$ performs ideally
  &
  depends on finer
  &
  $\Delta_{\textrm{LS}}$ outperforms $\Delta_{\textrm{oracle}}$\\
    &
  for any $\delta$
  &
  properties of $\pdf(x)$
  &
  for small $\delta < \delta_0$
  \\
  \hline
  Compressible
  & {\bf YES}
  & {\bf YES} or {\bf NO}
  & {\bf NO}\\
  \hline\hline
 &   &  {\em Proposition~\ref{prop:P0} (Section~\ref{sec: incompatible}):} & {\em Section~\ref{sec:GGD}:}   \\
 &   &  $\pdf_0(x) := 2 |x|/(x^2+1)^3$ & $\pdf_\tau(x) \propto \exp(-|x|^\tau)$ \\
 & & & $0<\tau<\infty$\\
 &   &  $\Delta_{\textrm{oracle}}$ performs just as  $\Delta_{\textrm{LS}}$&  Generalized Gaussian \\
 Examples & & & \\
 \cline{2-4}
  & \multicolumn{3}{c|}{\em Example~\ref{ex:sparseprior} (Section~\ref{sec:Discussion}):} \\
  &  \multicolumn{3}{c|}{$\pdf_{\tau,s}(x) \propto (1+|x|^\tau)^{-s/\tau}$}\\
 &   \multicolumn{3}{c|}{Generalized Pareto ($\tau=1$) / Student's $t$ ($\tau=2$)}\\
  \cdashline{2-4}
 &  \multicolumn{1}{c:}{ Case $1<s\leq 3$} &  \multicolumn{1}{c:}{Case $3<s<5$}  & Case $s>5$   \\
   &  \multicolumn{1}{c:}{} &  \multicolumn{1}{c:}{$\Delta_{\textrm{oracle}}$ outperforms $\Delta_{\textrm{LS}}$}  &\\
     & \multicolumn{1}{c:}{} & \multicolumn{1}{c:}{for small $\delta < \delta_0$}  &   \\
\hline
\end{tabular}
\end{table*}

\subsection{Are natural images compressible or incompressible ?}

Theorems~\ref{Th:4MomentIntro} and~\ref{th:InfiniteMomentInstOpt} provide easy to check conditions for (in)compressibility of a PDF $\pdf(x)$ based on its second of fourth moments. These rules of thumb are summarized in Table~\ref{table: summary}, providing an overview at a glance of the main results obtained in this paper.

We conclude this extended overview of the results with stylized application of these rules of thumb to wavelet and discrete cosine transform (DCT) coefficients of the natural images from the Berkeley database \cite{MartinFTM01}. Our results below provide an approximation theoretic perspective to the probabilistic modeling approaches in natural scene statistics community \cite{ruderman1994statistics,simoncelli2001natural,portilla2003image}.

\begin{figure*}\centering
\begin{tabular}{ccc}
{\includegraphics[width=0.3\hsize]{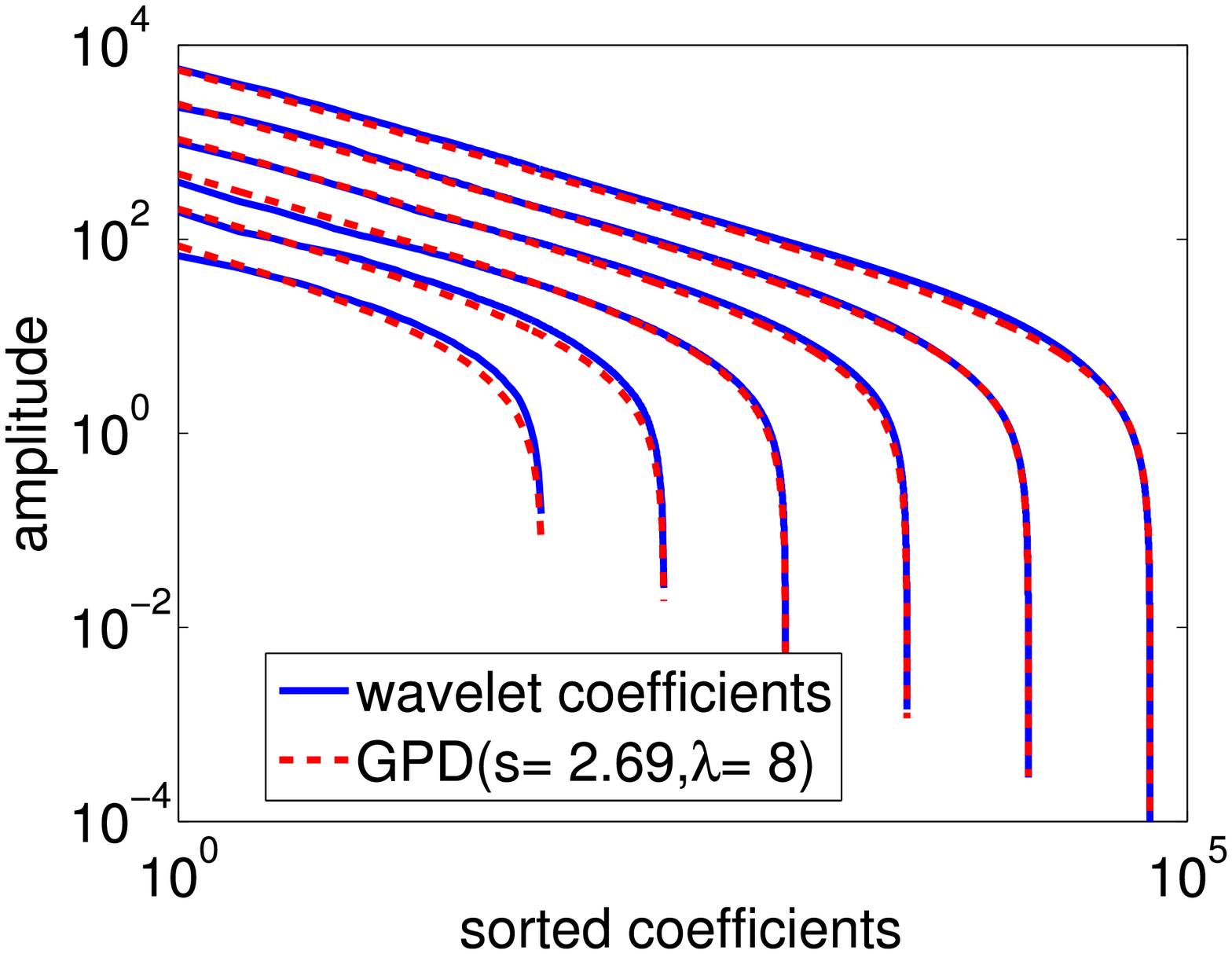}}&
{\includegraphics[width=0.3\hsize]{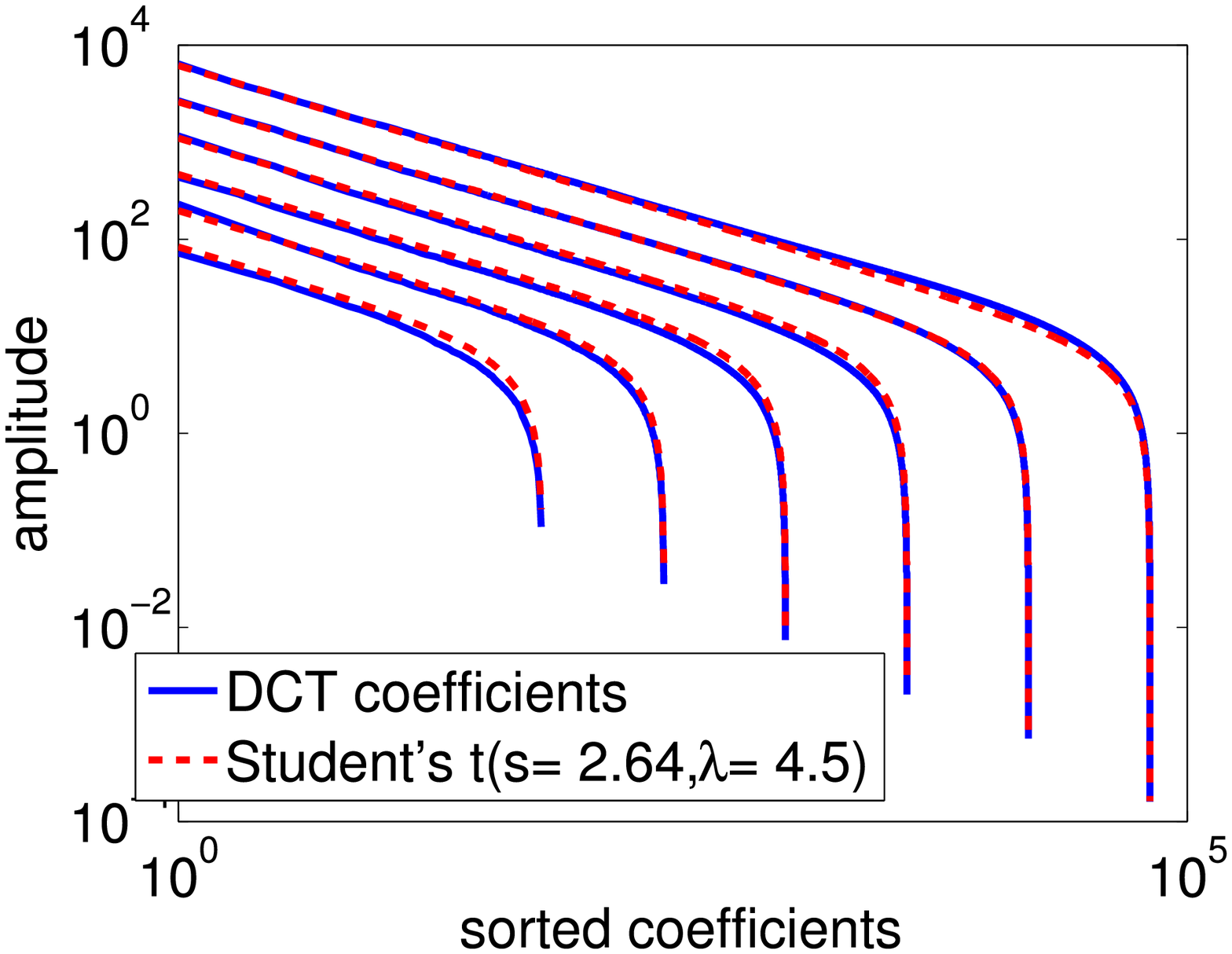}} &
{\includegraphics[width=0.32\hsize]{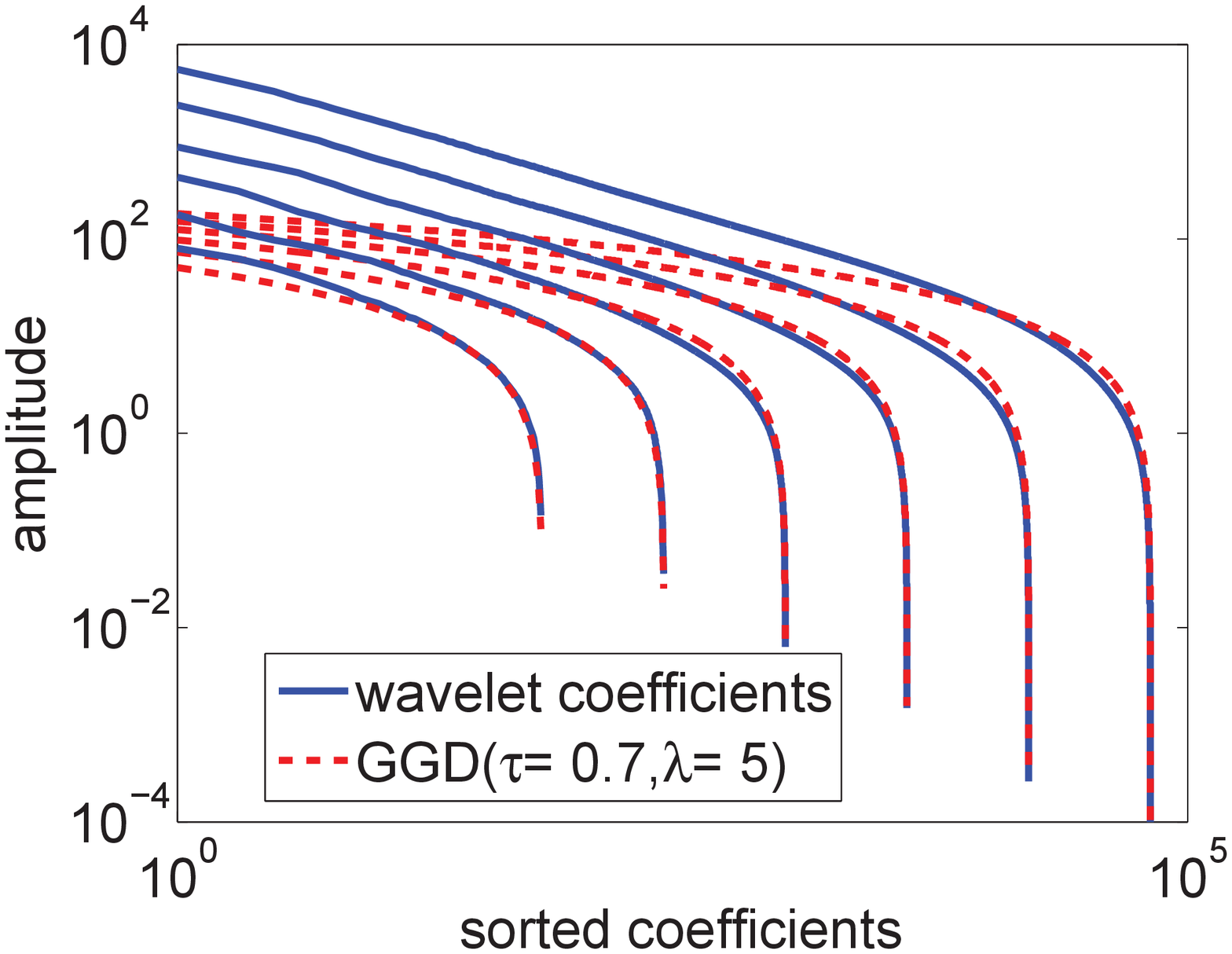}}\\
(a) Wavelet/GPD & (b) DCT/Student's $t$ distribution &(c) Wavelet/GGD
\end{tabular}
\caption{ \label{fig:waveDCT} Solid lines illustrate the Wavelet or DCT transform domain average order statistics of image patches from the Berkeley database \cite{MartinFTM01}. Dashed lines show the theoretical expected order statistics of the GPD, Student's $t$, and the GGD distributions with the indicated parameter values. The resolution of image patch sizes varies from left to right as $\{(8\times 8), (16\times 16),\ldots,(256\times 256)\}$, respectively.}
\end{figure*}

Figure~\ref{fig:waveDCT} illustrates, in log-log scale, the average of the magnitude ordered wavelet coefficients (Figures~\ref{fig:waveDCT}-(a)-(c)), and of the DCT coefficients (Figure~\ref{fig:waveDCT}-(b)).
They are obtained by randomly sampling $100$ image patches of varying sizes $N=2^j \times 2^j$ ($j=3,\ldots,8$), and taking their transforms (scaling filter for wavelets: Daubechies4).
For comparison, we also plot the expected order statistics (dashed lines), as described in \cite{IEEE_NIPS09_Cevher}, of the following distributions (cf Sections~\ref{sec:GGD} and~\ref{sec:Discussion})
\begin{itemize}
\item GPD: the scaled generalized Pareto distribution with density $\frac{1}{\lambda} \pdf_{\tau,s}(x/\lambda)$, $\tau=1$, with parameters
$s 
= 2.69$ and $\lambda = 8$ (Figure~\ref{fig:waveDCT}-(a));
\item Student's $t$: the scaled Student's $t$ distribution with density $\frac{1}{\lambda} \pdf_{\tau,s}(x/\lambda)$, $\tau=2$, with parameters
$s 
= 2.64$ and $\lambda = 4.5$ (Figure~\ref{fig:waveDCT}-(b));
\item GGD: the scaled generalized Gaussian distribution with density $\frac{1}{\lambda} \pdf_{\tau}(x/\lambda)$, with
$\tau = 0.7$ and $\lambda = 5$ (Figure~\ref{fig:waveDCT}-(c)).
\end{itemize}
The GGD parameters were obtained by approximating the histogram of the wavelet coefficients at $N=8\times 8$, as it is the common practice in the signal processing community \cite{Chang2000}.
The GPD and Student's $t$ parameters were tuned manually.

One should note that image transform coefficients are certainly not iid~\cite{Seeger2008}, for instance: nearby wavelets have correlated coefficients; wavelet coding schemes exploit well-known zero-trees indicating correlation across scales; the energy across wavelet scales often follows a power law decay.

The empirical goodness-of-fits in Figure \ref{fig:waveDCT} (a), (b)
seem to indicate that the distribution of the coefficients of natural
images, marginalized across all scales (in wavelets) or frequencies
(DCT) can be well approximated by a distribution of the type
$\pdf_{\tau,s}$ ({\em cf} Table~\ref{table: summary}) with
``compressibility parameter'' $s \approx 2.67 < 3$.  For this regime
the results of \cite{IEEE_NIPS09_Cevher} were inconclusive regarding
compressibility. However, from Table~\ref{table: summary} we see that such
a distribution satisfies
$\Expect X^2 = \infty$ ({\em cf} Example~\ref{ex:sparseprior} in
Section~\ref{sec:Discussion}), and therefore we are able to conclude that in the limit of very
high resolutions $N \to \infty$, such images are sufficiently
compressible to be acquired using compressive sampling with both {\em
  arbitrary good relative precision} and {\em arbitrary small
  undersampling} factor $\delta = m/N \ll 1$.

Considering the GGD with parameter $\tau = 0.7$, the results of Section~\ref{sec:GGD} ({\em cf} Figure~\ref{fig:GGcompressibility}) indicate that it is associated to a critical undersampling ratio $\delta_0(0.7) \approx 0.04$. Below this undersampling ratio, the oracle sparse decoder is outperformed by the least square decoder, which has the very poor expected relative error $1-\delta \geq 0.96$. Should the GGD be an accurate model for coefficients of natural images, this would imply that compressive sensing of natural images requires a number of measures at least $4\%$ of the target number of image pixels. However, while the generalized Gaussian approximation of the coefficients appear quite accurate at $N=8\times 8$, the empirical goodness-of-fits quickly deteriorate at higher resolution. For instance, the initial decay rate of the GGD coefficients varies with the dimension. Surprisingly, the GGD coefficients approximate the small coefficients (i.e., the histogram) rather well irrespective of the dimension. This phenomenon could be deceiving while predicting the compressibility of the images.


%

\section{Instance optimality, $\ell^r$-balls and compressibility in G-CS}
\label{sec: instance optimality}

Well-known results indicate that for certain matrices, $\sensing$, and for
certain types of sparse estimators of $\x$, such as the minimum
$\ell^1$ norm solution, $\Delta_1(\y)$,
an instance optimality property holds \cite{Cohen:2006aa}.
In the simplest case of noiseless observations, this reads:  the pair
$\{\sensing,\Delta\}$ is instance optimal to order $k$ in the $\ell^q$
norm with constant $C_k$ if for all $\x$:
\begin{equation}\label{eq:DefGeneralInstanceOptimality}
\|\Delta(\sensing\x)-\x\|_q \leq C_k \cdot \sigma_k(\x)_q
\end{equation}
where $\sigma_k(\x)_q$ is the error of best approximation of $\x$ with
$k$-sparse vectors, while $C_k$ is a constant which depends on $k$.
Various flavors of instance optimality are
possible~\cite{stablesigrecovery-CandesRombergTao-2006,Cohen:2006aa}.
We will initially focus on $\ell^1$ instance
optimality. For the $\ell^1$ estimator \eqref{eq:P1} it is known that instance optimality
in the $\ell^1$ norm (i.e. $q = 1$ in~\eqref{eq:DefGeneralInstanceOptimality}) is related to the following robust null space property. The matrix $\sensing$ satisfies the robust null space property of
order $k$ with constant $\eta \leq 1$ if:
\begin{equation}\label{eq:robustNSP}
||\z_\Omega||_1 < \eta ||\z_{\bar{\Omega}}||_1
\end{equation}
for all nonzero $\z$ belonging to the null space $\textrm{kernel} (\sensing) := \{\z, \sensing \z = 0\}$ and all index sets $\Omega$
of size $k$, where the notation $\z_\Omega$ stands for the vector matching $\z$ for indices in $\Omega$ and zero elsewhere.
It has further been shown~\cite{Davies:2009aa,Xu:2010} that the
robust null space property of order $k$ with constant $\eta_k$ is a necessary and sufficient condition for
$\ell^1$-instance optimality with the constant $C_k$
given by:
\begin{equation}\label{eq: instance optimality constant}
C_k = 2 \frac{(1+\eta_k)}{(1-\eta_k)}
\end{equation}

\begin{figure}[htbp]
\begin{center}
\includegraphics[width=.4\textwidth]{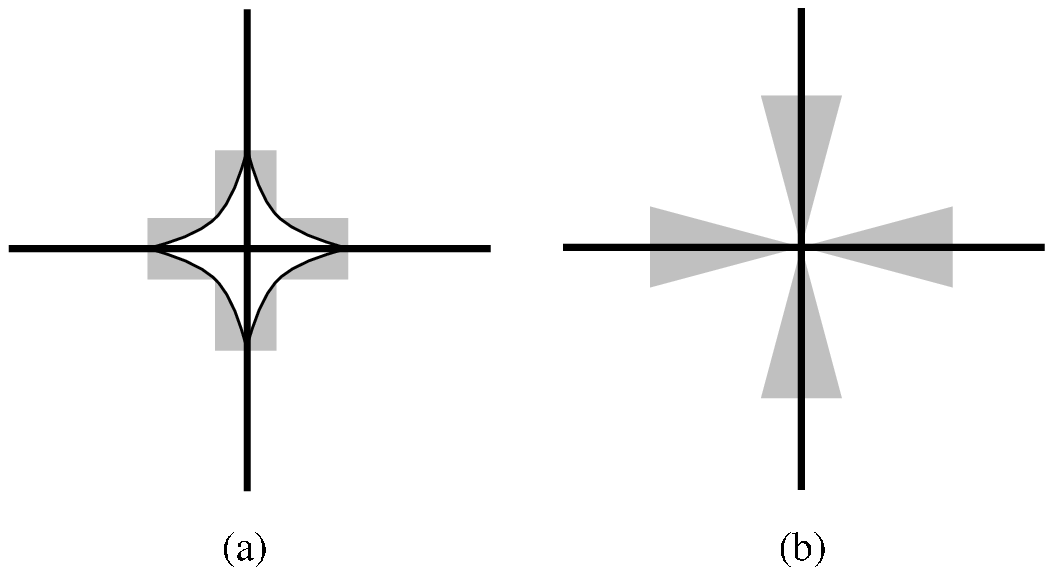}
\caption{\label{fig:lp_balls_pic}
(a) A cartoon view of an $\ell^r$ ball (white) and the weak
 $\ell^r$ ball of the same radius (grey); (b) A cartoon view of the
 notion of the compressible rays model.}
\end{center}
\end{figure}

Instance optimality is commonly considered as a strong property, since
it controls the {\em absolute} error in terms of the ``compressibility'' of
$\x$, expressed through $\sigma_k(\x)$. For instance optimality to be
meaningful we therefore require that $\sigma_k(\x)$ be small in some
sense. This idea has been encapsulated in a deterministic notion of
compressible vectors \cite{Cohen:2006aa}. Specifically suppose that
$\x$ lies in the $\ell^r$ ball of radius $R$ or the \emph{weak} $\ell^r$ ball of 
radius $R$ defined as:
\begin{equation}\label{eq:weak lp}
\|\x\|_{w\ell^r} := \sup_n \left\{|\x|_n^* \cdot n^{1/r}\right\} \leq R,
\end{equation}
with $|x|_n^*$ the $n$-th largest absolute value of elements of $\x$
(Figure~\ref{fig:lp_balls_pic}(a) illustrates the relationship between
the weak $\ell^r$ ball and the $\ell^r$ ball of the same radius). 
Then we can bound $\sigma_k(\x)_q$ for $q>r$, by
\begin{equation}
\sigma_k(\x)_q \leq R \left(\frac{r}{q-r}\right)^{1/q} k^{-(1/r-1/q)},
\end{equation}
therefore guaranteeing that the $k$-term approximation error is
vanishingly small for large enough $k$. 

Such models cannot be directly applied to the stochastic framework
since, as noted in \cite{Amini:2011}, iid realizations do not belong to
any weak $\ell^r$ ball. 
One obvious way to resolve this is to normalize the stochastic vector.
If $\Expect |X|^r = C < \infty$ then by the strong law of large
numbers, 
\begin{equation}\label{eq: normalized lp}
 \|\x_N\|_r^r/N \xrightarrow{a.s} C. 
\end{equation}
For example, such a signal model is considered in \cite{Donoho:2011}
for the G-CS problem,
where precise bounds on the worst-case asymptotic minimax mean-squared
reconstruction error are calculated for $\ell^1$ based decoders.

It can be tempting to assert that a vector drawn from a probability
distribution satisfying \eqref{eq: normalized lp} is ``compressible.'' Unfortunately, this is a
poor definition of a compressible distribution because finite dimensional $\ell^r$ balls also contain `flat'
vectors with entries of similar magnitude, that have very small
$k$-term approximation error \ldots only because the vectors are very
small themselves. 

\begin{figure}[htbp]
\begin{center}
\includegraphics[height=.5\textwidth,angle=90]{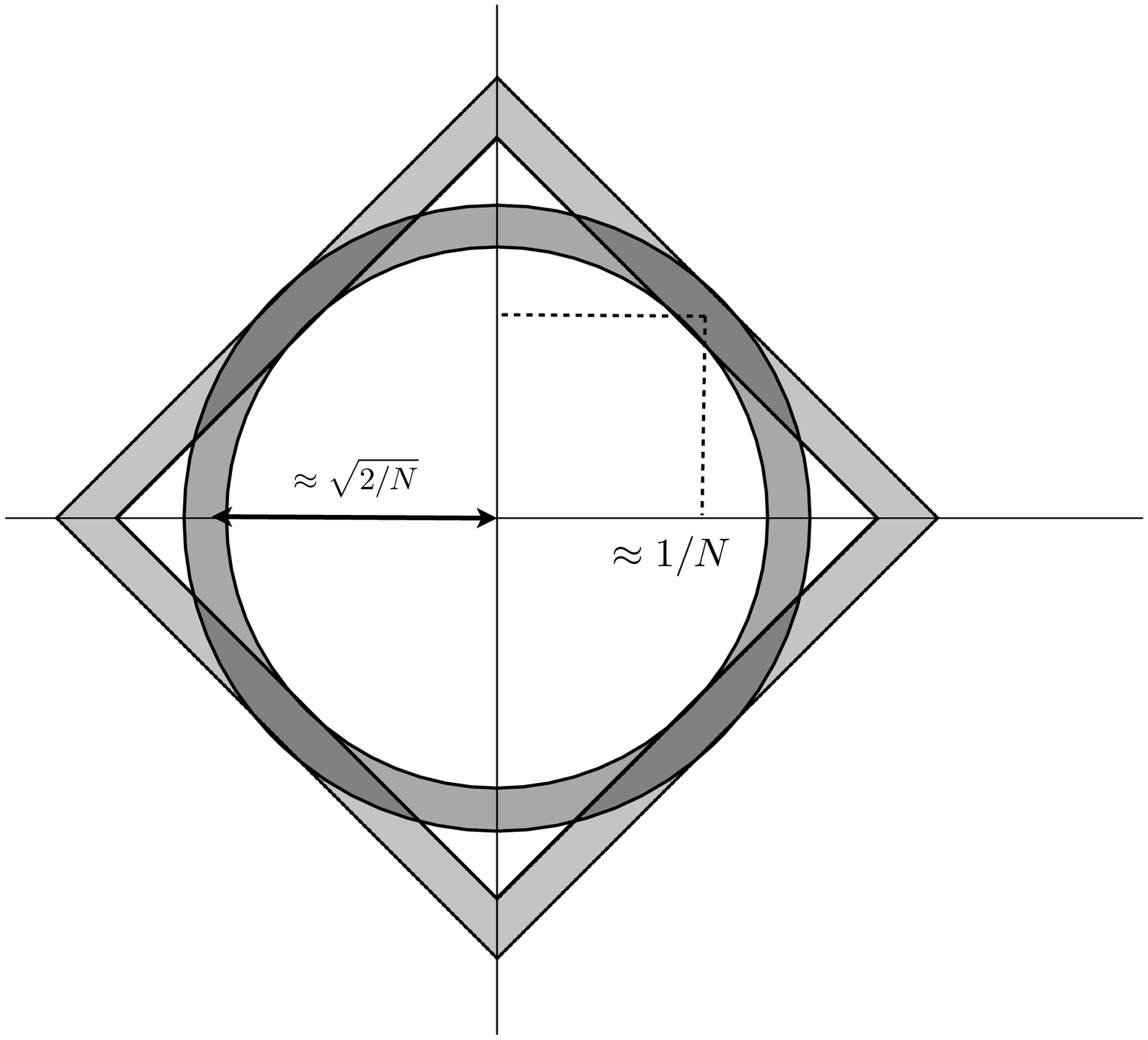}
\caption{A cartoon view of the $\ell^1$ and $\ell^2$ ``rings'' where vectors with iid Laplace-distributed entries concentrate. The radius of the $\ell^2$ ring is of the order of $\sqrt{2/N}$ while that of the $\ell^1$ ring is one, corresponding to vectors with flat entries $|\x|_n \approx 1/N$.}
\label{fig:ConcL1L2}
\end{center}
\end{figure}

For example, if $\x_N$ has entries drawn from the Laplace distribution then $\x_N/N$ will,  with high
probability, have an $\ell^1$-norm close to $1$. However the Laplace
distribution also has a finite second moment $\Expect X^2 =
2$, hence, with high probability $\x_N/N$ has $\ell^2$-norm  close to $\sqrt{2/N}$. This is not far from the $\ell^2$ norm of the largest flat
vectors that live in the unit $\ell^1$ ball, which have the form $|\x|_n
= 1/N$, $1 \leq n \leq N$. Hence a typical iid Laplace
distributed vector is a small and relatively flat vector. This is
illustrated on Figure~\ref{fig:ConcL1L2}. 

Instead of model \eqref{eq: normalized lp} we consider a more natural normalization of $\sigma_k(\x)_q$ with
respect to the size of the original vector $\x$ measured in the same
norm. This is the best $k$-term {\em relative} error
$\bar{\sigma}_k(\x)_q$ that we investigated in
Proposition~\ref{prop:RelativeError}. The class of vectors defined by 
$\bar{\sigma}_k(\x)_q < C$ for some $C$ does not have the shape
of an $\ell^r$ ball or weak $\ell^r$ ball. Instead it forms a set of
\emph{compressible `rays'} as depicted
in Figure~\ref{fig:lp_balls_pic} (b).

\subsection{Limits of G-CS  guarantees using instance optimality}

In terms of the relative best $k$-term approximation error, the
 instance optimality implies the following inequality:
\[
\frac{\|\Delta(\sensing \x)-\x\|}{\|\x\|}
\leq \min_k \left\{ C_k \cdot \bar{\bestapprox}_k(\x) \right\}
\]

Note that if we have the following inequality satisfied for the particular realization of $\x$
\[
\frac{\bestapprox_k(\x)}{\|\x\|} \geq C_k^{-1}, \forall k,
\]
then the only consequence of instance optimality is that
$\|\Delta(\sensing \x)-\x\| \leq \|\x\|$.
In other words, the performance
guarantee for the considered vector $\x$ is no better than for the trivial zero estimator:
$\Delta_{\textrm{trivial}}(\y) = 0$, for any $\y$.

This simple observation illustrates that one should  be careful in the
interpretation of instance optimality.  In particular, decoding
algorithms with instance optimality guarantees may not universally
perform better than other simple or more standard estimators.

To understand what this implies for specific distributions, consider the case
of $\ell^1$ decoding with a Gaussian encoder
$\Phi_N$. For this coder, decoder pair, $\{ \Phi_N, \Delta_1\}$, we know there is a strong phase transition associated
with the robust null space property~\eqref{eq:robustNSP} with $0<\eta<1$ (and hence the instance optimality
property with $1<C<\infty$) in terms of the undersampling factor $\mlevel := m/N$ and the
factor $\slevel : = k/m$ as $k,m,N \rightarrow \infty$
\cite{Xu:2010}. This is a generalization of the $\ell^1$ exact
recovery phase transition of Donoho and Tanner~\cite{Donoho:2009aa}
which corresponds to $\eta = 1$. We can therefore identify the
smallest instance optimality constant asymptotically possible as a
function of $\slevel$ and $\mlevel$ which we will term $C(\slevel,\mlevel)$.


To check whether instance optimality guarantees can beat the trivial zero
estimator $\Delta_{\textrm{trivial}}$ for a given undersampling ratio $\mlevel$, and a given
generative model $\pdf(x)$, we need to consider the product of
$\bar{\bestapprox}_k(\x)_1 \stackrel{a.s.}{\rightarrow} \Gfun{\pdf}{1}(\kappa)$
and $C(\frac{\kappa}{\mlevel},\mlevel)$. If
\begin{equation}
\Gfun{\pdf}{1}(\kappa)  > \frac{1}{C\left(\frac{\kappa}{\mlevel},\mlevel\right)},\quad \forall \kappa \in [0,\mlevel]
\end{equation}
then the instance optimality offers no guarantee to outperform
the trivial zero estimator.

In order to determine the actual strength of instance optimality we make the following
observations:
\begin{itemize}
\item $C(\frac{\kappa}{\mlevel},\mlevel) \geq 2$ for all $\kappa$
 and $\delta$;
\item $C(\frac{\kappa}{\mlevel},\mlevel) = \infty$ for all $\mlevel$ if $\kappa > \kappa_0 \approx 0.18$.
\end{itemize}
The first observation comes from minimising $C_k$ in \eqref{eq: instance
 optimality constant} with respect to $0 \leq \eta \leq 1$. The
 second observation stems from the fact that $\kappa_0 := \max_{\{\eta,\mlevel\}} \slevel_\eta(\mlevel)
 \approx 0.18$ \cite{Donoho:2009aa} (where $\slevel_\eta(\delta)$ is the strong threshold associated to the null space property with constant $\eta \leq 1$) therefore
 we have $\kappa = \mlevel \rho \leq \kappa_0 \approx 0.18$ for any finite $C$. 
 From these observations we obtain
:\\

{\em For distributions with PDF $\pdf(x)$ satisfying $\Gfun{\pdf}{1}(\kappa_0) \geq 1/2$, in high dimension $N$, 
instance optimality results for the decoder $\Delta_1$ with a Gaussian encoder can at best guarantee the performance (in the $\ell^{1}$ norm) 
of \ldots the trivial decoder $\Delta_{\textrm{trivial}}$.}\\


One might try to weaken the analysis by considering typical joint
behavior of $\Phi_N$ and $\x_N$. This
corresponds to the `weak' phase transitions
\cite{Donoho:2009aa,Xu:2010}. For this scenario there is a modified
$\ell^1$ instance optimality property \cite{Xu:2010}, however the constant
still satisfies $C(\frac{\kappa}{\mlevel},\mlevel) \geq
2$. Furthermore since $\kappa \leq \delta$ we can define an
undersampling ratio $\delta_0$ by $G_1[p](\delta_0) = 1/2$,
 such that weak instance optimality provides no guarantee that $\Delta_1$
will outperform the trivial decoder $\Delta_{\textrm{trivial}}$ in the
region $0< \delta \leq \delta_0$. More careful analysis will only
increase the size of this region.

\begin{example}[The Laplace distribution]
Suppose that $\x_N = (X_1,\ldots,X_N)$ has iid entries $X_n$ that follow the Laplace distribution with PDF $\pdf_1(x)$. Then for large $N$, as noted in Example~\ref{ex:Laplace}, the relative best $k$-term error is given by:
\[
\Gfun{\pdf_1}{1}(\kappa) = 1-\kappa \cdot \Big(1 + \ln 1/\kappa \Big)
\]
 Figure~\ref{fig:instanceOptimality_is_irrelevant} shows that
 unfortunately this function exceeds $1/2$ on the interval $\kappa \in
 [0,\ \kappa_0]$ indicating there are no non-trivial performance guarantees from
 instance optimality. Even exploiting weak instance optimality we can have
 no non-trivial guarantees below $\delta_0 \approx 0.18$.

\begin{figure}[htbp]
\begin{center}
\includegraphics[width=8cm]{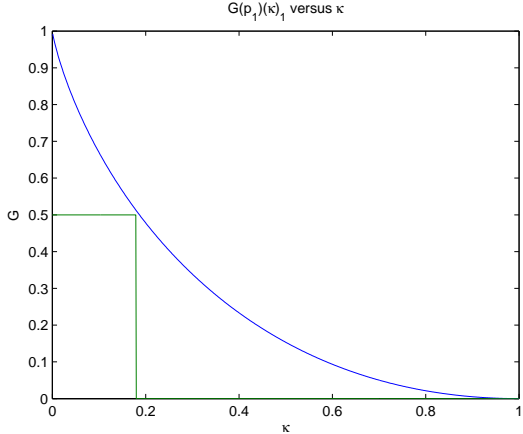}
\caption{The $\ell^1$-norm best $k$-term approximation relative error
 $\Gfun{\pdf_1}{1}(\kappa)$ as a function of $\kappa = k/N$ (top curve) along with a rectangular
 shaped function (bottom curve) that upper bounds $\inf_\delta C^{-1}(\kappa/\delta,\delta)$.}
\label{fig:instanceOptimality_is_irrelevant}
\end{center}
\end{figure}
\end{example}

\subsection{CS guarantees for random variables with unbounded second moment}

A more positive result (Theorem~\ref{th:InfiniteMomentInstOpt}) can be obtained showing that random
variables with infinite second moment, which are highly compressible (cf Proposition~\ref{prop:RelativeError}), are almost perfectly estimated by the $\ell^1$ decoder $\Delta_1$. In short, the result is based upon a variant of instance optimality: $\ell^2$
instance optimality in probability \cite{Cohen:2006aa} which can be shown to hold for
a large class of random matrices \cite{DeVore2009}. This can be
combined with the fact that when $\Expect
X^2 = \infty$, from Proposition~\ref{prop:RelativeError}, we have
$\Gfun{\pdf}{2}(\kappa) = 0$ for all $0 < \kappa \leq 1$ to give Theorem~\ref{th:InfiniteMomentInstOpt}. The proof is in the Appendix.


%
%
\begin{remark} A similar result can be derived
  based on $\ell^1$ instance optimality that shows that when $\Expect
 |X| = \infty$, then  the relative error in $\ell^1$ for the $\ell^1$
 decoder with a Gaussian encoder asymptotically goes to zero:
\begin{equation*}
\lim_{N \to \infty}\frac{\| \Delta_1(\sensing_N \x_N)-\x_N \|_1}{\|\x_N\|_1}
\stackrel{a.s.}{=} 0.
\end{equation*}
Whether other results hold for general $\ell^p$ decoders and relative $\ell^{p}$ error is not known.
\end{remark}

We can therefore conclude that a random variable with infinite
variance is not only compressible (in the sense of Proposition~\ref{prop:RelativeError}): it can also be accurately approximated from
undersampled measurements within a
compressive sensing scenario. In contrast, instance optimality provides no
guarantees of compressibility when the variance is finite and $\Gfun{\pdf}{1}(\kappa_0) \geq 1/2$.
At this juncture it is not clear where the blame for this result
lies. Is it in the strength of the instance optimality theory, or are
distributions with finite variance simply not able to generate sufficiently compressible vectors for sparse recovery to be successful at all?
We will explore this latter question further in subsequent
sections.



\section{G-CS performance of oracle sparse reconstruction {\em vs} least squares}
\label{sec:expectedperf}
Consider $\x$ an arbitrary vector in $\mathbb{R}^N$ and $\sensing$
be an $m \times N$ Gaussian encoder, and let $\y := \sensing \x$. Besides the trivial zero estimator $\Delta_{\textrm{trivial}}$ \eqref{eq:DeftrSolution} and the $\ell^1$ minimization estimator $\Delta_1$ \eqref{eq:P1}, the Least Squares (LS) estimator $\Delta_{\textrm{LS}}$ \eqref{eq:DefLSSolution} is a commonly used alternative.
Due to the Gaussianity of $\sensing$ and its independence from $\x$, it is well known that the resulting relative expected performance is
\begin{equation}
\label{eq:ExpectedPerfLS}
\frac{\Ephi \|\Delta_{\textrm{LS}}(\sensing \x) - \x\|_2^2} { \|\x\|_2^2}
= 1-\frac mN.
\end{equation}
Moreover, there is indeed a concentration around the expected value, as expressed by the inequality below: 
\begin{equation}
\label{eq:ConcentrationIneqLS}
(1-\epsilon) \left(1-\frac mN\right)
\leq
\frac{\|\Delta_{\textrm{LS}}(\sensing \x) - \x\|_2^2} { \|\x\|_2^2}
\leq
(1-\epsilon)^{-1} \left(1-\frac mN\right),
\end{equation}
for any $\epsilon > 0$ and $\x \in \R^N$, except with probability at most $2 \cdot e^{-(N-m) \epsilon^2/4} + 2 \cdot e^{-N \epsilon^2/4}$.

The result is independent of the vector $\x$, which
should be no surprise since the Gaussian distribution is
isotropic. The expected performance is directly governed by
the {\em undersampling factor}, i.e. the ratio between the number of measures $m$ and the dimension $N$ of the vector $\x$,
$\mlevel :=  m/N$.

In order to understand which statistical PDFs $\pdf(x)$ lead to ``compressible enough'' vectors $\x$, we wish to compare the performance of LS with that of estimators $\Delta$ that exploit the sparsity of $\x$ to estimate it.
Instead of choosing a particular estimator (such as $\Delta_1$), we consider the oracle sparse estimator $\Delta_{\textrm{oracle}}$ defined in~\eqref{eq:DefOSSolution}, which is likely to upper bound the performance of most sparsity based estimators. While in practice $\x$ must be estimated from $\y = \sensing \x$, the oracle is given a precious side information: the index set $\Lambda$ associated to the $k$ largest components in $\x$, where $k < m$. Given this information, the oracle computes
\begin{equation*}
\Delta_{\textrm{oracle}}(\y,\Lambda) := \argmin_{\textrm{support}(\x) = \Lambda} \|\y-\sensing \x\|_2^2 = \sensing_{\Lambda_k}^+ \y,
\end{equation*}
where, since $k<m$, the pseudo-inverse is
$\sensing_{\Lambda}^+ =
(\sensing_{\Lambda}^T\sensing_{\Lambda})^{-1}
\sensing_{\Lambda}^T$.
Unlike LS,  the expected performance of the oracle estimators drastically depend on the shape of the best $k$-term approximation relative error of $\x$. Denoting $\x_I$ the vector whose entries match those of $\x$ on an index set $I$ and are zero elsewhere, and $\bar{I}$ the complement of an index set, we have the following result.
\begin{theorem}[Expected performance of oracle sparse estimation]
\label{th:ExpectedPerfOracle}
Let $\x \in \R^N$ be an arbitrary vector, $\sensing$ be an $m \times N$ random Gaussian matrix, and $\y := \sensing \x$. Let $\Lambda$ be an index set of size $k < m-1$, either deterministic, or random but {\em statistically independent} from $\sensing$. We have
\begin{align}
\frac{\Ephi \|\Delta_{\textrm{oracle}}(\sensing\x,\Lambda)-\x\|_2^2}{ \|\x\|_2^2}
&=
\frac{1}{1-\frac{k}{m-1}} \times
\frac{\|\x_{\bar{\Lambda}}\|_2^2}{\|\x\|_2^2}\notag\\
\label{eq:OracleKTermError}
&\geq
\frac{1}{1-\frac{k}{m-1}} \times
\frac{\sigma_k(\x)_2^2}{\|\x\|_2^2}.
\end{align}
If $\Lambda$ is chosen to be the $k$ largest components
of $\x$, then the last inequality is an equality.
Moreover, we can characterize the concentration around the expected value as
\begin{align}
\label{eq:OracleKTermConc}
1+\frac{k(1-\epsilon)^3}{m-k+1}
&\leq
\frac{\|\Delta_{\textrm{oracle}}(\sensing\x,\Lambda)-\x\|_2^2}{ \|\x_{\bar{\Lambda}}\|_2^2}
\leq
1+\frac{k(1-\epsilon)^{-3}}{m-k+1}\\
\intertext{except with probability at most}
\label{eq:OracleKTermConc1}
8 \cdot & e^{-\min(k,m-k+1) \cdot c_l(\epsilon)/2},
\intertext{where}
\ c_l(\epsilon) & := -\ln (1-\epsilon) - \epsilon \geq \epsilon^2/2.
\end{align}
\end{theorem}

\begin{remark}
Note that this result assumes that $\Lambda$ is {\em statistically independent} from $\sensing$.
Interestingly, for practical decoders such as the $\ell^1$ decoder, $\Delta_1$, the selected $\Lambda$
might not satisfy this assumption, unless the decoder successfully identifies the support of the largest components of $\x$.
\end{remark}

\subsection{Compromise between approximation and conditioning}
We observe that the expected performance of both $\Delta_{LS}$ and $\Delta_{\textrm{oracle}}$
is essentially governed by the quantities $\mlevel = m/N$ and $\slevel = k/m$, which are
reminiscent of the parameters in the phase transition diagrams of Donoho and
Tanner~\cite{Donoho:2009aa}. However, while in the work of Donoho and
Tanner the quantity $\slevel$ parameterizes a {\em model} on the
vector $\x_N$, which is assumed to be $\slevel\mlevel N$-sparse, here
$\slevel$ rather indicates the order of $k$-term approximation of
$\x_N$ that is {\em chosen} in the oracle estimator.
In a sense, it is more related to a stopping criterion that one would use in a greedy algorithm. The
quantity that actually models $\x_N$ is the function
$\Gfun{\pdf}{2}$, provided that $\x_N \in \R^N$ has iid entries $X_n$ with PDF $\pdf(x)$ and finite second moment $\Expect X^2 < \infty$. Indeed, combining Proposition~\ref{prop:RelativeError} and Theorem~\ref{th:ExpectedPerfOracle} we obtain:

\begin{theorem}
\label{th:ASCVOracle}
Let $\x_{N}$ be iid with respect to $\pdf(x)$ as in Proposition~\ref{prop:RelativeError}. Assume that $\Expect X^2 < \infty$. Let $\phi_{i,j}$, $i,j \in \mathbb{N}$ be iid Gaussian variables $\mathcal{N}(0,1)$. Consider two sequences $k_N, m_N$ of integers and assume that
\begin{equation}
\lim_{N \to \infty} k_N/m_N = \slevel\quad \mbox{and}\quad \lim_{N \to \infty} m_N/N = \mlevel.
\end{equation}
Define
the $m_N \times N$ Gaussian encoder
$\sensing_N = \left[\phi_{ij}/\sqrt{m_N} \right]_{1 \leq i \leq m_N, 1\leq j \leq N}$. Let $\Lambda_N$ be the index of the $k_N$ largest magnitude coordinates of $\x_N$. We have the almost sure convergence
\begin{eqnarray}
\label{eq:ASCVOracle}
\lim_{N \to \infty}
\frac{\|\Delta_{\textrm{oracle}}(\sensing_N \x_{N},\Lambda_N)-\x_N\|_2^2}{ \|\x_N\|_2^2}
&\stackrel{a.s.}{=}& \frac{\Gfun{\pdf}{2}(\slevel\mlevel)}{1-\slevel};\\
\label{eq:ASCVLS}
\lim_{N \to \infty}
\frac{\|\Delta_{\textrm{LS}}(\sensing_N\x_{N})-\x_N\|_2^2}{ \|\x_N\|_2^2}
&\stackrel{a.s.}{=}& 1-\delta.
\end{eqnarray}
\end{theorem}

For a given undersampling ratio $\mlevel = m/N$, the asymptotic expected performance of the oracle
therefore depends on the relative number of components that are kept $\slevel =k/m$, and we observe the same tradeoff as discussed in Section~\ref{sec: instance optimality}:
\begin{itemize}
\item For large $k$, close to the number of measures $m$ ($\slevel$ close to one),
the ill-conditioning of the pseudo-inverse matrix $\sensing_\Lambda$ (associated to the factor
$1/(1-\slevel)$) adversely impacts the expected performance;
\item For smaller $k$, the pseudo-inversion of this matrix
is better conditioned, but the $k$-term approximation error governed by
$\Gfun{\pdf}{2}(\slevel \mlevel)$ is increased.
\end{itemize}
Overall, for some intermediate size  $k  \approx \slevel^\star m$ of the oracle support set $\Lambda_k$,
the best tradeoff between good approximation and good conditioning is achieved, leading at best to the asymptotic expected performance
\begin{equation}
\label{eq:DefBestSparsityConditioningTradeoff}
\Hfun{\pdf}(\mlevel)
:=  \inf_{\slevel \in (0,1)} \frac{\Gfun{\pdf}{2}(\slevel \mlevel)}{1-\slevel}.
\end{equation}


\section{A comparison of least squares and oracle sparse methods}
\label{sec:ComparisonL2L0}
The question that we will now investigate is how the expected performance of oracle sparse methods compares to that of least squares, i.e., how large is $\Hfun{\pdf}(\mlevel)$ compared to $1-\mlevel$?
We are particularly interested in understanding how they compare for
small $\mlevel$.  Indeed, large $\mlevel$ values are associated with
scenarii that are quite irrelevant to, for example, compressive sensing since the
projection $\sensing \x$ cannot significantly compress the dimension
of $\x$. Moreover,  it is in the regime where $\mlevel$ is small that
the expected performance of least squares is very poor, and we would
like to understand for which PDFs $\pdf$ sparse
approximation is an inappropriate tool.
The answer will of course depend on the PDF $\pdf$ through the
function $\Gfun{\pdf}{}(\cdot)$.  To characterize this we will say
that a PDF $\pdf$ is \emph{incompressible at a subsampling rate of
  $\mlevel$} if
\[ \Hfun{\pdf}(\mlevel) > 1-\mlevel.
\]
In practice, there is often a minimal undersampling
rate, $\mlevel_0$, such that for $\mlevel \in (0,\mlevel_0)$ least
squares estimation dominates the oracle sparse estimator. Specifically
we will show below that PDFs $\pdf(x)$ with a finite fourth moment
$\Expect X^4 < \infty$, such as generalized Gaussians, always
have some minimal undersampling rate $\mlevel_0 \in (0,1)$ below which
they are incompressible. As a result, unless we perform at least
$m \geq \mlevel_0  N$ random Gaussian measurement of an
associated $\x_N$, it is not worth relying on sparse methods for
reconstruction since least squares can do as good a job.

When the fourth moment of the
distribution is infinite, one might hope that the converse is true, i.e.
that no such minimal
undersampling rate $\mlevel_0$  exists. However, this is not the
case. We will show that there is a PDF $\pdf_0$, with infinite
fourth moment and finite second moment, such that
\[
\Hfun{\pdf_0}(\mlevel) = 1-\mlevel,\qquad \forall \mlevel \in (0,1).
\]
Up to a scaling factor, this PDF is associated to the symmetric PDF
\begin{equation}
\label{eq:DefP0}
\pdf_0(x) := \frac{2|x|}{(x^2+1)^3}
\end{equation}
and illustrates that least squares can be competitive with oracle
sparse reconstruction even when the fourth moment is infinite.

\subsection{Distributions incompatible with extreme undersampling}
\label{sec: incompatible}
In this section we show that when a PDF $\pdf(x)$ has a finite fourth
moment, $\Expect X^4 < \infty$, then it will generate vectors which are not sufficiently compressible to be compatible with compressive sensing at high level of undersampling. We begin by
showing that the comparison of $\Hfun{\pdf}(\mlevel)$ to
$1-\mlevel$ is related to that of $\Gfun{\pdf}{2}(\kappa)$ with
$(1-\sqrt{\kappa})^2$.

\begin{lemma}
\label{le:Hbound}
Consider a function $G(\kappa)$ defined on $(0,1)$ and define
\begin{equation}
H(\mlevel) := \inf_{\slevel \in (0,1)} \frac{G(\mlevel\slevel)}{1-\slevel}.
\end{equation}
\begin{enumerate}
\item If  $G(\mlevel^2) \leq (1-\mlevel)^2$, \\ then $H(\mlevel)\leq1-\mlevel$.
\item If $G(\kappa) \leq (1-\sqrt{\kappa})^2$ for all $\kappa \in (0,\sqrt{\mlevel_0})$, \\ then
$H(\mlevel) \leq 1-\mlevel$ for all $\mlevel \in (0, \mlevel_0)$.
\item If $G(\kappa) \geq (1-\sqrt{\kappa})^2$ for all $\kappa \in (0,\mlevel_0)$, \\ then $H(\mlevel) \geq 1-\mlevel$ for all $\mlevel \in (0, \mlevel_0)$.
\end{enumerate}
\end{lemma}

Lemma~\ref{le:Hbound} allows us to deal directly with
$\Gfun{\pdf}{2}(\kappa)$ instead of
$\Hfun{\pdf}(\mlevel)$. Furthermore the $(1-\sqrt{\kappa})^2$
term can be related to the fourth moment of the distribution (see Lemma~\ref{le:G4Moment} in the Appendix) giving
the following result, which implies Theorem~\ref{Th:4MomentIntro}:

\begin{theorem}
\label{th:4Moment}
If $\Expect_{\pdf(x)} X^4 < \infty$, then there exists a minimum
undersampling $\mlevel_0 = \mlevel_0[\pdf] > 0$ such that
for $\mlevel < \mlevel_0$,
\begin{equation}
\Hfun{\pdf}(\mlevel) \geq 1-\mlevel, \forall\ \mlevel \in (0,\mlevel_0).
\end{equation}
and the performance of the oracle $k$-sparse estimation as
described in Theorem~\ref{th:ASCVOracle} is asymptotically
almost surely worse than that of least squares estimation as $N
\rightarrow \infty$.
\end{theorem}

Roughly speaking, if $\pdf(x)$ has a finite fourth moment, then in the
regime where the relative number of measurement is (too) small we obtain a
better reconstruction with least squares than with the oracle sparse
reconstruction!

Note that this is rather strong, since the oracle is allowed to know
not only the support of the $k$ largest components of the unknown
vector, but also the best choice of $k$ to balance approximation error
against numerical conditioning.
A striking example is the case of generalized Gaussian distributions
discussed below.

One might also hope that, reciprocally, having an infinite fourth
moment would suffice for a distribution to be compatible with compressed sensing at extreme levels of undersampling. The following result disproves this hope.
\begin{prop}
\label{prop:P0}
With the PDF $\pdf_0(x)$ defined in~\eqref{eq:DefP0}, we have
\begin{equation}
\Hfun{\pdf_0}(\mlevel) = 1-\mlevel, \forall\ \mlevel \in (0,1).
\end{equation}
\end{prop}

On reflection this should not be that surprising. The
PDF $\pdf_0(x)$ has no probability mass at $x = 0$ and
resembles a smoothed Bernoulli distribution with heavy tails.

\subsection{Worked example: the generalized Gaussian distributions}
\label{sec:GGD}
Theorem~\ref{th:4Moment} applies in particular whenever $\x_N$ is
drawn from a generalized Gaussian distribution,
\begin{equation}
\label{eq:GGD}
\pdf_\tau(x) \propto \exp\left(-c |x|^\tau\right),
\end{equation}
where $0 < \tau < \infty$.  The shape parameter, $\tau$
controls how heavy or light the tails of the distribution are. When
$\tau = 2$ the distribution reduces to the standard Gaussian, while for $\tau <2$
it gives a family of heavy tailed distributions with positive
kurtosis. When $\tau = 1$ we have the Laplace distribution and for
$\tau \leq 1$ it is often considered that the distribution is in some way ``sparsity-promoting''.
However, the generalized Gaussian always has a finite fourth
moment for all $\tau>0$. Thus Theorem~\ref{th:4Moment} informs us that for a given parameter $\tau$
there is always a critical undersampling value below which the
generalized Gaussian is incompressible.

\begin{figure}[htbp]
\begin{center}
\includegraphics[width=8cm]{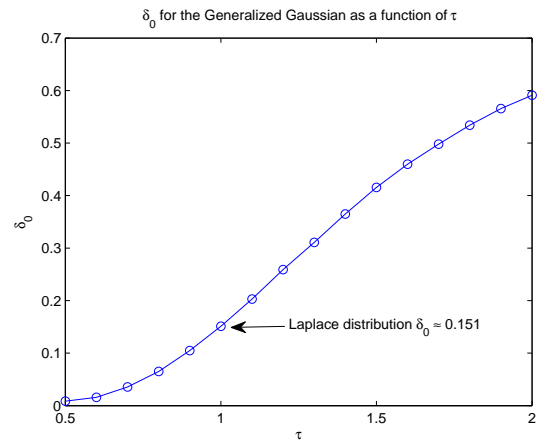}
\caption{A plot of the critical subsampling rate, $\mlevel_0$ below which the
  generalized Gaussian distribution is incompressible as a
  function of the shape parameter, $\tau$.\label{fig:GGcompressibility}
}
\end{center}
\end{figure}

While Theorem~\ref{th:4Moment} indicates the existence of a critical
$\mlevel_0$ it does not provide us with a useful bound. Fortunately,
although in general we are unable to derive explicit expressions for
$\Gfun{\pdf}{}(\cdot)$ and $\Hfun{\pdf}(\mlevel)$ (with the
exceptions of $\tau = 1,2$ - see Appendix~\ref{app:Laplace}), the
generalized Gaussian has a closed form expression for its cdf in
terms of the incomplete gamma function.
\[
\cdf(x) = \frac{1}{2}+ \sgn(x) \frac{\gamma\left(1/\tau,c|x|^\tau
\right)}{2 \Gamma(1/\tau)}
\]
where $\Gamma(\cdot)$ and $\gamma(\cdot,\cdot)$
 are respectively the gamma function and the lower incomplete gamma function.
 We are therefore able to
numerically compute the value of $\mlevel_0$ as a function of
$\tau$ with relative ease. This is shown in
Figure~\ref{fig:GGcompressibility}.
We see that, unsurprisingly, when $\tau$ is around $2$ there is little to
be gained even with an oracle sparse estimator over standard least
squares estimation. When $\tau = 1$ (Laplace distribution) the value
of $\mlevel_0 \approx 0.15$, indicating that when subsampling by a
factor of roughly 7 the least squares estimator will be
superior. At this level of undersampling the relative
error is a very poor: $0.85$, that is a performance of $0.7$dB in terms of
traditional Signal to Distortion Ratio (SDR).

The critical undersampling value steadily drops as $\tau$ tends
towards zero and
the distribution becomes increasingly leptokurtic. Thus data
distributed according to the generalized
Gaussian for small $\tau \ll 1$ may still be a reasonable candidate for
compressive sensing distributions as long as the undersampling rate is kept
significantly above the associated $\mlevel_0$.

\subsection{Expected Relative Error for the Laplace distribution}

We conclude this section by examining in more detail the performance
of the estimators for Laplace distributed data at various
undersampling values. We have already seen from
Figure~\ref{fig:GGcompressibility} that the oracle performance is poor
when subsampling by roughly a factor of 7. What about more modest
subsampling factors? Figure~\ref{fig:LaplaceRelError} plots the
relative error as a
function of undersampling rate, $\mlevel$. The horizontal lines indicate SDR
values of $3$dB, $10$dB and $20$dB. Thus for the oracle estimator to achieve 10dB
the undersampling rate must be greater than 0.7, while to achieve a
performance level of $20$dB, something that might reasonably be
expected in many sensing applications, we can hardly afford any
subsampling at all since this requires $\mlevel > 0.9$.

At this point we should remind the reader that these performance
results are for the comparison between the {\em oracle} sparse estimator and
linear least squares. For practically implementable reconstruction algorithms we
would expect that the critical undersampling rate at which least
squares wins would be significantly higher. Indeed, as shown in
Figure~\ref{fig:LaplaceRelError}, this is what is empirically observed
for the average performance of the $\ell^1$ estimator \eqref{eq:P1}
applied to Laplace distributed data. This curve was calculated at
various values of $\mlevel$ by averaging the relative error of 5000 $\ell^1$
reconstructions of independent Laplace distributed realizations of
$\x_N$ with $N = 256$. In particular note that the
$\ell^1$ estimator only outperforms least squares for undersampling $\delta$
above approximately 0.65!



\section{Concluding discussion}\label{sec:Discussion}
As we have just seen, Generalized Gaussian distributions are
incompressible at low subsampling rates because their fourth moment is
always finite. This confirms the results of Cevher obtained with a
different approach~\cite{IEEE_NIPS09_Cevher}, but may come as a
surprise: for $0 < \tau \leq 1$ the minimum $\ell^\tau$ norm solution
to $\y = \sensing \x$, which is also the MAP estimator under the
Generalized Gaussian prior, is known to be a good estimator of $\x_0$
when $\y = \sensing \x_0$ and $\x_0$ is compressible
\cite{Davies:2009aa}.  This highlights the need to distinguish between
an estimator and its MAP interpretation. In contrast, we
describe below a family of PDFs $\pdf_{\tau,s}$ which, for certain values of the parameters $\tau,s$, combines:
\begin{itemize}
\item  superior asymptotic almost sure performance of oracle sparse estimation over least squares reconstruction $\Delta_\textrm{oracle}$, even
in the largely undersampled scenarios $\mlevel \to 0$;
\item connections between oracle sparse estimation and MAP estimation.
\end{itemize}

\begin{example}
\label{ex:sparseprior}
For $0<\tau <\infty,\ 1<s<\infty$ consider the probability density function
\begin{equation}
\pdf_{\tau,s}(x) \propto (1+|x|^\tau)^{-s/\tau}.
\end{equation}
\begin{enumerate}
\item {\bf When $1 < s\leq 3$, the distribution is compressible.}\\ 
Since $\Expect_{\pdf_{\tau,s}} X^2 = \infty$, Theorem~\ref{th:InfiniteMomentInstOpt} is applicable: the $\ell^1$ decoder with a Gaussian encoder has ideal asymptotic performance, even at arbitrary small undersampling $\delta = m/N$;
\item {\bf When $3<s<5$, the distribution remains somewhat compressible.}\\
On the one hand $\Expect_{\pdf_{\tau,s}} X^2 < \infty$, on the other hand $\Expect_{\pdf_{\tau,s}} X^4 = \infty$. \\
A detailed examination of the $\Gfun{\pdf_{\tau,s}}{1}$ function shows that 
there exists a relative number of measures $\mlevel_0(\tau,s) > 0$ such that in the low measurement regime $\mlevel < \mlevel_0$,
the asymptotic almost sure performance of oracle of $k$-sparse estimation, as described
in Theorem~\ref{th:ASCVOracle}, with the best choice of $k$, is {\em better} than that of least squares estimation:
\begin{equation}
\Hfun{\pdf_{\tau,s}}(\mlevel) < 1-\mlevel, \forall \mlevel \in (0,\mlevel_0).
\end{equation}
\item {\bf When $s > 5$,  the distribution is incompressible.}\\
Since $\Expect_{\pdf_{\tau,s}} X^4 < \infty$, Theorem~\ref{Th:4MomentIntro} is applicable: with a Gaussian encoder, there is an  undersampling ratio $\delta_0$ such that whenever $\delta < \delta_0$, the asymptotic almost sure performance of oracle sparse estimation is worse than that of least-squares estimation;
\end{enumerate}
\end{example}
Comparing Proposition~\ref{prop:P0} with the above
Example~\ref{ex:sparseprior}, one observes that both the PDF
$\pdf_0(x)$ (Equation~\eqref{eq:DefP0}) and the PDFs $\pdf_{\tau,s}$,
$3<s<5$ satisfy $\Expect_{\pdf_{\tau,s}} X^2 < \infty$ and
$\Expect_{\pdf_{\tau,s}} X^4 = \infty$. Yet, while $\pdf_0$ is
essentially incompressible, the PDFs $\pdf_{\tau,s}$ in this range are
compressible. This indicates that, for distributions with finite
second moment and infinite fourth moment, compressibility depends not
only on the tail of the distribution {\em but also on their mass
  around zero}. However the precise dependency is currently unclear.

For $\tau = 2$, the PDF $\pdf_{2,s}$ is a Student-t distribution. For $\tau = 1$, it is called a generalized Pareto
distribution. These have been considered in~\cite{IEEE_NIPS09_Cevher,Baraniuk:2010aa} as examples of
``compressible'' distributions, with the added condition that $s \leq 2$. Such a restriction results from the use of $\ell^2-\ell^1$ instance optimality in~\cite{IEEE_NIPS09_Cevher,Baraniuk:2010aa}, which implies that sufficient compressibility conditions can only be satisfied when $\Expect_\pdf |X| = \infty$. Here instead we exploit $\ell^2-\ell^2$ instance optimality {\em in probability}, making it possible to obtain compressibility when $\Expect X^2 = \infty$.  In other words,~\cite{IEEE_NIPS09_Cevher,Baraniuk:2010aa} provides {\em sufficient} conditions on a PDF $\pdf$ to check its compressibility, but is inconclusive in characterizing their incompressibility.

The family of PDFs, $\pdf_{\tau ,s}$ in the range $0< \tau \leq1$, can also be linked with a
sparsity-inducing MAP estimate.
Specifically for an observation $\y = \sensing \x$ of a given vector
$\x \in \R^N$,
one can define the MAP estimate under the probabilistic
{\em model} where all entries of $\x$ are considered as
iid distributed according to $\pdf_{\tau,s}$:
\[
\Delta_{\textrm{MAP}}(\y) := \arg\max_{\x | \sensing \x = \y} \prod_{n=1}^N \pdf_{\tau,s}(x_n) =
\argmin_{\x | \sensing \x = \y} \sum_{n=1}^N f_\tau(|x_n|).
\]
where for $t \in \R^+$ we define $f_\tau(t) :=  \log (1+ t^\tau ) = a_{\tau,s} -b_{\tau,s} \log\pdf_{\tau,s}(|t|)$.
One can check that the function $f_\tau$ is associated to an admissible $f$-norm as described in~\cite{gribonval04:highlysparsICA04,gribonval07:_highl}: $f(0) = 0$, $f(t)$ is non-decreasing, $f(t)/t$ is
non-increasing (in addition, we have $f(t) \sim_{t \to 0}  t^\tau$).
Observing that the MAP estimate is a ``minimum $f$-norm'' solution to
the linear problem $\y = \sensing \x$, we can conclude that whenever
$\x$ is a ``sufficiently (exact) sparse'' vector, we have in
fact~\cite{gribonval04:highlysparsICA04,gribonval07:_highl}
$\Delta_{\textrm{MAP}}(\sensing\x) = \x$, and $\Delta_{\textrm{MAP}}(\sensing\x) = \Delta_1(\sensing \x)$ is also the minimum $\ell^1$ norm solution to $\y = \sensing \x$, which can in turn be ``interpreted'' as
the MAP estimate under the iid Laplace model.  However, unlike the
Laplace interpretation of $\ell^1$ minimization, here
Example~\ref{ex:sparseprior} indicates that such densities are better
aligned to sparse reconstruction techniques. Thus the MAP estimate
interpretation here may be more valid.

It would be interesting to determine
whether the MAP estimator $\Delta_{\textrm{MAP}}(\sensing\x)$ for such
distributions is in some way close to optimal (i.e. close to the
minimum mean squared error solution for $\x$). This would give such
estimators a degree of legitimacy from a Bayesian
perspective. However, we have \emph{not} shown that the estimator
$\Delta_{\textrm{MAP}}(\sensing\x)$ provides a good estimate for data
that is distributed according to $\pdf_{\tau ,s}$ since, if $\x$ is a
large dimensional typical instance with
entries drawn iid from the PDF $\pdf_{\tau,s}(x)$,  it is typically
{\em not} exactly sparse, hence the uniqueness results
of~\cite{gribonval04:highlysparsICA04,gribonval07:_highl} do not
directly apply.  One would need to resort to a more detailed
robustness analysis in the spirit of~\cite{gribonval06:simpletest} to
get more precise statements relating
$\Delta_{\textrm{MAP}}(\sensing\x)$ to $\x$.


%
\appendix
\subsection{Proof of Proposition~\ref{prop:RelativeError}} 
\label{app:Wald}
To prove Proposition~\ref{prop:RelativeError} we will rely on the following theorem \cite{Bruss:1991aa}[Theorem 2.2].
\begin{theorem}
\label{th:Wald}
Suppose that $\cdf_Y$ is a continuous and strictly increasing cumulative density function on $[a,b]$ where $0 \leq a < b \leq \infty$, with $\cdf_Y(a)=0$, $\cdf_Y(b)=1$. For $\sigma \in (0,\mu)$ where $\mu = \int_a^b y d\cdf_Y(y)$, let $\tau \in (a,b)$ be defined by the equation $\sigma = \int_a^\tau y d\cdf_Y(y)$. Let $s_1,s_2,\ldots$ be a sequence such that $\lim_{N\to \infty} s_N/N = \sigma$, and let $Y_1,Y_2 \ldots$ be iid random variables with cumulative density function $\cdf_Y$. Let $Y_{1,N} \leq \ldots \leq Y_{N,N}$ be the increasing order statistics of $Y_1,\ldots, Y_N$ and let $L_N = L(N,s_n)$ be defined as
$L(N,s_n) := 0$ if $Y_{1,N} > s_N$, otherwise:
\begin{align}
\label{eq:DefLnWald}
L(N,s_n) := \max \left\{\ell \leq N, Y_{1,N} + \ldots + Y_{\ell,N} \leq s_N \right\};
\end{align}
Then
\begin{align}
\label{eq:WaldResult}
&\lim_{N \to \infty} \frac{Y_{L_N,N}}{N} 
\stackrel{a.s.}{=} \tau,\\
&\lim_{N \to \infty} \frac{L_N}{N} 
\stackrel{a.s.}{=} \cdf_Y(\tau),\\
&\lim_{N \to \infty} \frac{\mathbb{E}(L_N)}{N} 
= \cdf_Y(\tau).
\end{align}
\end{theorem}

\begin{IEEEproof}[Proof of Proposition~\ref{prop:RelativeError}]
We begin by the case where $\Expect |X|^q < \infty$. We consider random variables $X_n$ drawn according the PDF $\pdf(x)$, and we define the iid non-negative random variables $Y_n = |X_n|^q$. They have the cumulative density function $\cdf_Y(y) = \mathbb{P}(Y\leq y) = \mathbb{P}(|X|\leq y^{1/q}) = \bar{\cdf}(y^{1/q})$, and we have $\mu = \Expect Y = \Expect |X|^q = \int_0^\infty |x|^q d\bar{\cdf}(x) \in (0,\infty)$. We define $\mathbf{x}_N = (X_n)_{n=1}^N$, and we consider a sequence $k_N$ such that $\lim_{N \to \infty} k_N/N = \kappa \in (0,1)$. By the assumptions on $\cdf_Y$ there is a unique $\tau_0 \in (0,\infty)$ such that $\kappa = 1-\cdf_Y(\tau_0)$, and we will prove that
\begin{eqnarray}
\label{eq:WaldConsequenceLimInf}
\liminf_{N \to \infty} \frac{\sigma_{k_N}(\mathbf{x}_N)_q^q}{N \mu}
&\stackrel{a.s.}{\geq}& \frac{\int_0^{\tau_0} y d\cdf_Y(y) }{\mu},\\
\label{eq:WaldConsequenceLimSup}
\limsup_{N \to \infty} \frac{\sigma_{k_N}(\mathbf{x}_N)_q^q}{N \mu}
&\stackrel{a.s.}{\leq}& \frac{\int_0^{\tau_0} y d\cdf_Y(y) }{\mu}.
\end{eqnarray}
The proof of the two bounds is identical, hence we only detail the first one. Fix $0<\epsilon<\tau_0$ and define $\tau = \tau(\epsilon) := \tau_0-\epsilon$, $\sigma = \sigma(\epsilon) := \int_0^\tau y d\cdf_Y(y)$, and $s_N = N\sigma$. Defining $L_N$ as in~\eqref{eq:DefLnWald}, we can apply Theorem~\ref{th:Wald} and obtain $\lim_{N \to \infty} \frac{L_N}{N} \stackrel{a.s.}{=} \cdf_Y(\tau)$. Since $\lim_{N \to \infty} \frac{k_N}{N} = 1-\cdf_Y(\tau_0)$, it follows that
\[
\lim_{N \to \infty} \frac{N-k_N}{L_N} \stackrel{a.s.}{=} \frac{\cdf_Y(\tau_0)}{\cdf_Y(\tau)} >1
\]
where we used the fact that $\cdf_Y$ is strictly increasing and $\tau<\tau_0$. In other words, almost surely, we have $N-k_N > L_N$ for all large enough $N$. Now remember that by definition
\[
L_N
= \max\left\{ \ell \leq N, \sigma_{N-\ell}(\mathbf{x}_N)_q^q \leq N \sigma \right\}.
\]
As a result, almost surely, for all large enough $N$, we have
\[
\sigma_{k_N}(\mathbf{x}_N)_q^q = \sigma_{N-(N-k_N)}(\mathbf{x}_N)_q^q > N \sigma.
\]
Now, by the strong law of large number, we also have
\[
\lim_{N \to \infty} \frac{\|\mathbf{x}_N\|_q^q}{N\mu} \stackrel{a.s.}{=} 1,
\]
hence we obtain
\[
\liminf_{N \to \infty} \frac{\sigma_{k_N}(\mathbf{x}_N)_q^q}{\|\mathbf{x}_N\|_q^q}
\stackrel{a.s.}{\geq} \frac{\sigma}{\mu} =
\frac{\int_0^{\tau_0-\epsilon} y d\cdf_Y(y) }{\mu}.
\]
Since this holds for any $\epsilon>0$ and $\cdf_Y$ is continuous, this implies~\eqref{eq:WaldConsequenceLimInf}.  The other bound~\eqref{eq:WaldConsequenceLimSup} is obtained similarly. Since the two match, we get
\[
\lim_{N \to \infty} \frac{\sigma_{k_N}(\mathbf{x}_N)_q^q}{\|\mathbf{x}_N\|_q^q}
\stackrel{a.s.}{=}
\frac{\int_0^{\tau_0} y d\cdf_Y(y) }{\mu}
=
\frac{\int_0^{\tau_0} y d\cdf_Y(y) }{\int_0^\infty y d\cdf_Y(y) }.
\]
Since $\kappa = 1-\cdf_Y(\tau_0) = 1-\bar{\cdf}(\tau_0^{1/q})$ we have $\tau_0 = \left[\bar{\cdf}^{-1}(1-\kappa)\right]^q$.
Since $\cdf_Y(y) = \bar{\cdf}(y^{1/q})$ we have $d \cdf_Y(y) = \frac{1}{q} y^{1/q-1} \bar{\pdf}(y^{1/q}) dy$. As a result
\begin{eqnarray*}
\frac{\int_0^{\tau_0} y d\cdf_Y(y) }{\int_0^\infty y d\cdf_Y(y) }
&=&
\frac{\int_0^{ \left[\bar{\cdf}^{-1}(1-\kappa)\right]^q} y^{1/q} \bar{\pdf}(y^{1/q}) dy }{\int_0^\infty y^{1/q} \bar{\pdf}(y^{1/q}) dy }\\
&\stackrel{(a)}{=}&
\frac{\int_0^{ \bar{\cdf}^{-1}(1-\kappa)} x \bar{\pdf}(x) x^{q-1} dx }{\int_0^\infty x \bar{\pdf}(x) x^{q-1} dx }\\
&=&
\frac{\int_0^{ \bar{\cdf}^{-1}(1-\kappa)} x^q \bar{\pdf}(x) dx }{\int_0^\infty x^q \bar{\pdf}(x) dx }
\end{eqnarray*}
where in (a) we used the change of variable $y = x^q$, $x = y^{1/q}$, $dy = q x^{q-1} dx$.
We have proved the result for $0<\kappa<1$, and we let the reader check that minor modifications yield the results for $\kappa=0$ and $\kappa=1$.

Now we consider the case $\Expect |X|^q = + \infty$. The idea is to use a ``saturated'' version $\tilde{X}$ of the random variable $X$, such that $\Expect |\tilde{X}|^q < \infty$, so as to use the results proven just  above.

One can easily build a family of smooth saturation functions $f_\eta: [0\ +\infty) \to [0\ 2\eta)$, $0 < \eta < \infty$ with $f_\eta(t) = t$, for $t \in [0,\eta]$, $f_\eta(t)\leq t$, for $t>\eta$, and two additional properties:
\begin{enumerate}
\item each function $t \mapsto f_\eta(t)$ is bijective from $[0,\infty)$ onto $[0,2\eta)$, with $f'_\eta(t) > 0$ for all $t$;
\item each function $t \mapsto f_\eta(t)/t$ is monotonically decreasing;
\end{enumerate}
Denoting $f_\eta(\x) := (f_\eta(x_i))_{i=1}^N$, by~\cite[Theorem 5]{gribonval07:_highl}, the first two properties ensure that
for all $1 \leq k \leq N$, $\x \in \R^N$, $0< \eta,q < \infty$ we have
\begin{equation}
\label{eq:BoundRelErrorWithSaturation}
\frac{\sigma_k(\x)^q}{\|\x\|_q^q}
\leq
\frac{\sigma_k(f_\eta(\x))^q}{\|f_\eta(\x)\|_q^q}.
\end{equation}
Consider a fixed $\eta$ and the sequence of ``saturated'' random variables $\tilde{X}_i = f_\eta(|X_i|)$. They are iid with $\Expect |\tilde{X}|^q < \infty$. Moreover, the first property of $f_\eta$ above ensures that their cdf $t \mapsto \bar{\cdf}_\eta(t) := \Prob(f_\eta(|X|)\leq t)$ is continuous and strictly increasing on $[0\ 2\eta]$, with $\bar{F}_\eta(0) = 0$ and $\bar{F}_\eta(\infty) = 1$. Hence, by the first part of Proposition~1 just proven  above, we have
\begin{align}
\label{eq:ASCVRelErrorWithSaturation}
\lim_{N \to \infty} \frac{\sigma_{k_N}(f_\eta(\x_N))^q}{\|f_\eta(\x_N)\|_q^q}
&\stackrel{a.s.}{=}
\Gfun{\bar{\pdf}_\eta}{q}(\kappa) =
\frac{\int_0^{\bar{\cdf}_\eta^{-1}(1-\kappa)} x^q
 \bar{\pdf}_\eta(x) dx}{\int_0^\infty x^q \bar{\pdf}_\eta(x) dx}\notag\\
& \leq
 \frac{|\bar{\cdf}_\eta^{-1}(1-\kappa)|^q}{\Expect |f_\eta(X)|^q}.
 \end{align}
Since $f_\eta(t) \leq t$ for all $t$, we have $\bar{F}_\eta(t) = \Prob(f_\eta(|X|) \leq t) \geq \Prob(|X| \leq t) = \bar{\cdf}(t)$ for all $t$, hence $\bar{\cdf}_\eta^{-1}(1-\kappa) \leq \bar{\cdf}^{-1}(1-\kappa)$. Moreover, since $f_\eta(t) = t$ for $0 \leq t \leq \eta$, we obtain $\Expect |f_\eta(X)|^q \geq \int_0^\eta x^q \bar{\cdf}(x) dx$. Combining~~\eqref{eq:BoundRelErrorWithSaturation} and~\eqref{eq:ASCVRelErrorWithSaturation} with the above observations we obtain for any $0<\eta<\infty$
\[
\limsup_{N \to \infty} \frac{\sigma_{k_N}(\x_N)^q}{\|\x_N\|_q^q}
\leq
\lim_{N \to \infty} \frac{\sigma_{k_N}(f_\eta(\x_N))^q}{\|f_\eta(\x_N)\|_q^q}
\stackrel{a.s.}{\leq} \frac{|\bar{\cdf}^{-1}(1-\kappa)|^q}{\int_0^\eta x^q \bar{\pdf}(x) dx}.
\]
Since $\Expect |X|^q = \int_0^\infty x^q \bar{\pdf}(x)dx = \infty$, the infimum over $\eta$ of the right hand side is zero.
\end{IEEEproof}

\begin{remark}
To further characterize the typical asymptotic behaviour of the relative error when
$\Expect_{\pdf}(|X|^q) = \infty$ and $k_N/N \to 0$ appears to require a
more detailed characterization of the probability density function,
such as decay bounds on the tails of the distribution.
\end{remark}



\subsection{Proof of Theorem~\ref{th:InfiniteMomentInstOpt}}

The proof is based upon the following version of 
 \cite[Theorem~$5.1$]{DeVore2009}: 

\begin{theorem}[DeVore \emph{et al.} \cite{DeVore2009}] \label{th:DeVore2009} 
Let $\sensing(\omega) \in \R^{m \times N}$ be a random matrix whose
entries are iid and drawn from ${\cal N}(0,1/m)$. There are some
absolute constants $C_0,\ldots, C_6$, and $C_7$ depending on $C_1,
\ldots, C_6$ such that, given  any $k \leq C_0
m/\log(N/m)$ then 
\begin{equation}
\label{eq: instoptinprob}
\| \x-\Delta_1(\sensing(\omega)\x)\|_2 \leq C_7
\bestapprox_k(\x)_2,
\end{equation} 
with probability exceeding 
\[1-C_1 e^{-C_2 m} -e^{-C_3 \sqrt{Nm}} 
- C_4 e^{-C_5 m}-2m e^{-\frac{\sqrt{m}}{C_6 \log(N/m)}}.
\]
\end{theorem}
In this version of the theorem we have specialized to the case where
the random matrices are Gaussian distributed. We have also removed the
rather peculiar requirement in the original version that $N \geq [\ln 6]^2 m$ as careful scrutiny of the proofs
(in particular the proof of Theorem 3.5 \cite{DeVore2009}) indicates that the effect of this term can be absorbed into the constant $C_3$ as
long as $m/N \leq [\tfrac{2}{\ln 6}]^2 \approx 1.2$, which is trivially satisfied.

We now proceed to prove
Theorem~\ref{th:InfiniteMomentInstOpt}. 
By assumption the undersampling ratio $\delta = \lim_{N
  \rightarrow \infty} \tfrac{m_N}{N} > 0$, therefore there exists a $0
  < \kappa < 1$ such that 
\begin{equation}
\notag
\delta > C_0 \kappa \log \frac{1}{\delta}.
\end{equation}
Now choosing a sequence $k_N/N \rightarrow \kappa$ we have, for large enough $N$,
\begin{align}
\notag
m_N &\geq C_0 k_N \log(N/m_N).
\end{align}
Hence, applying Theorem~\ref{th:DeVore2009}, for all $N$ large enough, there exist a set
$\Omega_N(\x_N,k_N)$ with
\begin{equation}
\Prob (\Omega_N^c(\x_N,k_N)) \leq C_8 m e^{-C_9 \sqrt{m}}
\end{equation}
such that~\eqref{eq: instoptinprob} holds for all $\sensing_N(\omega) \in \Omega(\x_N,k_N)$, i.e.,
\begin{equation}
 \frac{\| \x_N-\Delta_1(\Phi_N(\omega)\x_N)\|_2}{\|\x_N\|_2} \leq C_7 \bar{\bestapprox}_{k_N}(\x_N)_2.
\end{equation}
A union bound argument similar to the one used in the proof of Theorem~\ref{th:ASCVOracle} (see Appendix~\ref{app:PfThASCVOracle}) gives:
\begin{align}
\limsup_{N \rightarrow \infty} 
 \frac{\| \x_N-\Delta_1(\Phi_N\x_N)\|_2}{\|\x_N\|_2} &\stackrel{a.s.}{\leq} \limsup_{N
 \rightarrow \infty} C_7 \bar{\bestapprox}_{k_N}(\x_N)_2 \notag\\
&\stackrel{a.s.}= C_7  G_2[p](\kappa)
= 0.
\end{align}


\subsection{Proof of Theorem~\ref{th:ExpectedPerfOracle}}
We will need concentration bounds for several distributions. For the Chi-square distribution with $n$ degrees of freedom  $\chi^2_n$, we will use the following standard result (see, e.g., \cite[Proposition 2.2]{Barvinok:2005aa}, and the intermediate estimates in the proof of~\cite[Corollary 2.3]{Barvinok:2005aa}):
\begin{prop}
\label{prop:GaussConc}
Let $X \in \R^n$ a standard Gaussian random variable. Then, for any $0<\epsilon <1$
\begin{align}
\Prob\big(\|X\|_2^2\geq n(1-\epsilon)^{-1}\big)
& \leq e^{-n \cdot c_u(\epsilon)/2}\\
\Prob\big(\|X\|_2^2 \leq n (1-\epsilon)\big)
\leq e^{-n \cdot c_l(\epsilon)/2}
\end{align}
with
\begin{align}
c_u(\epsilon)
&:=
\frac{\epsilon}{1-\epsilon} + \ln(1-\epsilon)\\
c_l(\epsilon)
&:=
-\ln(1-\epsilon)-\epsilon.
\end{align}
\end{prop}
Note that
\begin{equation}
\epsilon^2/2 \leq c_l(\epsilon)  \leq c_u(\epsilon),\quad 0<\epsilon<1.
\end{equation}
Its corollary, which provides concentration for projections of random variables from the unit sphere, will also be useful. The statement is obtained by adjusting~\cite[Lemma 3.2]{Barvinok:2005aa} and~\cite[Corollary 3.4]{Barvinok:2005aa} keeping the sharper estimate from above.
\begin{corollary}
\label{cor:KDimSphere}
Let $X$ be a random vector uniformly distributed on the unit sphere in $\R^n$, and let $X_L$ be its orthogonal projection on a $k$-dimensional subspace $L$  (alternatively, let $X$ be an arbitrary random vector and $L$ be a random $k$-dimensional subspace uniformly distributed on the Grassmannian manifold). For any $0<\epsilon<1$ we have
\begin{align}
\Prob\big(\sqrt{\frac{n}{k}}\|X_L\|_2 \geq \|X\|_2(1-\epsilon)^{-1}\big)
 \leq e^{-k \cdot c_u(\epsilon)/2} + e^{-n \cdot c_l(\epsilon)/2},\\
\Prob\big(\sqrt{\frac{n}{k}}\|X_L\|_2 \leq \|X\|_2 (1-\epsilon)\big)
 \leq e^{-k \cdot c_l(\epsilon)/2} + e^{-n \cdot c_u(\epsilon)/2}.
\end{align}
\end{corollary}
The above result directly implies the concentration inequality~\eqref{eq:ConcentrationIneqLS} for the LS estimator mentioned in Section~\ref{sec:expectedperf}.
We will also need a result about Wishart matrices. The Wishart distribution~\cite{Muirhead:2008aa} $\mathcal{W}_\ell(n,\Sigma)$ is the distribution of $\ell \times \ell$ matrices $A = Z^T Z$ where $Z$ is an $n \times \ell$ matrix whose columns have the normal distribution $\mathcal{N}(0,\Sigma)$.
\begin{theorem}[{\cite{Muirhead:2008aa} [Theorem 3.2.12 and consequence, p. 97-98]}]
\label{th:InvWishart}
If $A$ is $\mathcal{W}_\ell(n, \Sigma)$ where $n-\ell+1>0$, and if $Z \in \R^\ell$ is a random vector distributed independently of $A$ and with $P(Z=0) = 0$, then the ratio $Z^T \Sigma^{-1}Z / Z^T A^{-1}Z$ follows a Chi-square distribution with $n-\ell+1$ degrees of freedom $\chi^2_{n-\ell+1}$, and is independent of $Z$. Moreover, if $n-\ell-1>0$ then
\begin{equation}
\label{eq:MuirheadExpectInvWishart}
\Expect A^{-1} = \Sigma^{-1} \cdot (n-\ell-1)^{-1}.
\end{equation}
\end{theorem}
Finally, for convenience we formalize below some useful but simple facts that we let the reader check.
\begin{lemma}
\label{le:SphericalIndependence}
Let $\mathbf{A}$ and $\mathbf{B}$ be two independent $m \times k$ and $m\times \ell$ random Gaussian matrices with iid entries $\mathcal{N}(0,1/m)$, and let $x \in \R^\ell$ be a random vector independent from $\mathbf{B}$. Consider a singular value decomposition (SVD)  $\mathbf{A} = U \Sigma V$ and let $u_\ell$ be the columns of $U$. Define $w := \mathbf{B} x/\|\mathbf{B}x\|_2 \in \R^m$, $w_1 := (\langle u_\ell,w\rangle)_{\ell=1}^k \in \R^k$, $w_2 := w_1/\|w_1\|_2 \in \R^k$ and $w_3 := V^T w_2 \in \R^k$. We have
\begin{enumerate}
\item\label{it:w} $w$ is uniformly distributed on the sphere in $\R^m$, and statistically independent from $\mathbf{A}$;
\item\label{it:w1} the distribution of $w_1$ is rotationally invariant in $\R^k$, and it is statistically independent from $\mathbf{A}$;
\item\label{it:w2} $w_2$ is uniformly distributed on the sphere in $\R^k$, and statistically independent from $\mathbf{A}$;
\item\label{it:w3} $w_3$ is uniformly distributed on the sphere in $\R^k$, and statistically independent from $\mathbf{A}$.
\end{enumerate}
\end{lemma}

We can now start the proof of Theorem~\ref{th:ExpectedPerfOracle}.
For any index set $J$, we denote $\x_J$ the
vector which is zero out of $J$. For matrices, the notation
$\sensing_J$ indicates the sub-matrix of $\sensing$ made of the
columns indexed by $J$. The notation $\bar{J}$ stands for the
complement of the set $J$.
For any index set $\Lambda$ associated to linearly independent
columns of $\sensing_\Lambda$ we can write
$\y= \sensing_\Lambda \x_\Lambda + \sensing_{\bar{\Lambda}} \x_{\bar{\Lambda}}$ hence
\begin{eqnarray}
\Delta_{\textrm{oracle}}(\y,\Lambda)  &:=&
\sensing_\Lambda^+ \y = \x_\Lambda + \sensing_\Lambda^+
\sensing_{\bar{\Lambda}} \x_{\bar{\Lambda}}\notag\\
\label{eq:DefOracleKTermError}
\|\Delta_{\textrm{oracle}}(\y,\Lambda) -\x\|_2^2
&=&
\|\sensing_\Lambda^+  \sensing_{\bar{\Lambda}}
\x_{\bar{\Lambda}}\|_2^2 + \|\x_{\bar{\Lambda}}\|_2^2
\end{eqnarray}
The last equality comes from the fact that the restriction of
$(\Delta_{\textrm{oracle}}(\y,\Lambda) -\x)$ to the indices in $\Lambda$ is  $\sensing_\Lambda^+
\sensing_{\bar{\Lambda}} \x_{\bar{\Lambda}}$, while its restriction to $\bar{\Lambda}$ is $\x_{\bar{\Lambda}}$.
Denoting
\begin{equation}
w := \frac{\sensing_{\bar{\Lambda}} \x_{\bar{\Lambda}}}{\|\sensing_{\bar{\Lambda}} \x_{\bar{\Lambda}}\|_2} \in \R^m
\end{equation}
we obtain the relation
\begin{equation}
\label{eq:DeterministicDecomp}
\frac{\|\Delta_{\textrm{oracle}}(\y,\Lambda) -\x\|_2^2}{\|\x_{\bar{\Lambda}}\|_2^2}
=
\underbrace{\|\sensing_\Lambda^+  w\|_2^2}_{A}
\times
\underbrace{\frac{\|\sensing_{\bar{\Lambda}}  \x_{\bar{\Lambda}}\|_2^2}{\|\x_{\bar{\Lambda}}\|_2^2}}_{B}
 +1.
\end{equation}
From the singular value decomposition
\[
\sensing_\Lambda = U_m \cdot
\left[
\begin{array}{c}
\Sigma_k\\
0_{(m-k) \times k}
\end{array}
\right]
\cdot V_k,
\]
where $U_m$ is an $m \times m$ unitary matrix with columns $u_\ell$, and $V_k$ is a $k \times k$ unitary matrix, we deduce that $\sensing_\Lambda^+ = V_k^T [\Sigma_k^{-1},0_{k \times (m-k) }] U_m^T$ and
\begin{equation}
\label{eq:DecomPhiPinv1}
\|\sensing_\Lambda^+  w\|_2^2
=
 \|[\Sigma_k^{-1} 0_{k \times (m-k) }] U_m^T w\|_2^2
= \sum_{\ell=1}^k \sigma_\ell^{-2} |\langle u_\ell, w \rangle|^2.
\end{equation}

Since $\sensing_{\bar{\Lambda}}$ and $\x_{\bar{\Lambda}}$ are statistically independent, the random vector $\sensing_{\bar{\Lambda}} \x_{\bar{\Lambda}} \in \R^m$ is Gaussian with zero-mean and covariance $m^{-1} \cdot \|\x_{\bar{\Lambda}}\|_2^2 \cdot \mathbf{Id}_m$. Therefore,
\begin{equation}
\label{eq:ExpectFactorA}
\Expect \left\{\|\sensing_{\bar{\Lambda}}\x_{\bar{\Lambda}}\|_2^2/\| \x_{\bar{\Lambda}}\|_2^2\right\} = 1
\end{equation}
and by Proposition~\ref{prop:GaussConc}, for any $0<\epsilon_0<1$
\begin{equation}
\label{eq:ConcFactorA}
\Prob\left(1-\epsilon_0 \leq \frac{\|\sensing_{\bar{\Lambda}}\x_{\bar{\Lambda}}\|_2^2}{\| \x_{\bar{\Lambda}}\|_2^2}
 \leq (1-\epsilon_0)^{-1} \right)
\geq
1-2 \cdot e^{-m \cdot c_l(\epsilon_0)/2}.
\end{equation}
Moreover, by Lemma~\ref{le:SphericalIndependence}-item~\ref{it:w1}, the random variables $\langle u_\ell,w\rangle$, $1 \leq \ell \leq k$ are identically distributed and independent from the random singular values $\sigma_\ell$. Therefore,
\begin{align*}
\Expect \|\sensing_\Lambda^+  w\|_2^2
=&
\Expect\left\{
\sum_{\ell = 1}^k
\sigma_\ell^{-2}
\right\}
\times \Expect \{|\langle u,w\rangle|^2\}\\
=& \Expect \left\{\text{Trace} (\sensing_\Lambda^T \sensing_\Lambda)^{-1} \right\} \times \frac 1m.
\end{align*}
The matrix $\sensing_\Lambda^T \sensing_\Lambda$ is $\mathcal{W}_k(m, \frac1m \mathbf{Id}_k)$ hence, by Theorem~\ref{th:InvWishart}, when $m-k-1>0$ we have
\begin{equation}
\label{eq:ExpectFactorB}
\Expect \|\sensing_\Lambda^+  w\|_2^2
=
\frac{\operatorname{Trace}(m\mathbf{Id}_k)}{(m-k-1) m} = \frac{k}{m-k-1}.
\end{equation}
Now, considering $w_1 := (\langle u_\ell,w\rangle)_{\ell=1}^k \in \R^k$, $w_2 := w_1/\|w_1\|_2$ and $w_3 := V_k^T w_2$, we obtain
\begin{align*}
\|\sensing^+ w\|_2^2
&= \|\Sigma_k^{-1} w_1\|_2^2 = \|w_1\|_2^2 \times \|\Sigma_k^{-1} w_2\|_2^2\\
&=\|w_1\|_2^2 \times \|\Sigma_k^{-1} V_k w_3\|_2^2\\
&=\|w_1\|_2^2 \times w_3^T (\sensing_\Lambda^T \sensing_\Lambda)^{-1} w_3 = m \|w_1\|_2^2/R(w_3),
\end{align*}
where $R(w_3) := m \|w_3\|_2^2/w_3^T (\sensing_\Lambda^T \sensing_\Lambda)^{-1} w_3 = w_3^T (m^{-1} \mathbf{Id}_k)^{-1} w_3/w_3^T (\sensing_\Lambda^T \sensing_\Lambda)^{-1} w_3$.
By Lemma~\ref{le:SphericalIndependence}-item~\ref{it:w3}, $w_3$ is statistically independent from $\sensing_\Lambda$. As a result, by Theorem~\ref{th:InvWishart}, the random variable $R(w_3)$ follows a Chi-square distribution with $m-k+1$ degrees of freedom $\chi^2_{m-k+1}$, and by Proposition~\ref{prop:GaussConc}, for any $0<\epsilon_1<1$,
\begin{align}
\label{eq:RatioBound}
\Prob&\Big(1-\epsilon_1 \leq R(w_3)^{-1} \cdot (m-k+1) \leq (1-\epsilon_1)^{-1} \Big)\notag\\
& \geq 1- 2 e^{-(m-k+1) \cdot c_l(\epsilon_1)/2}.\\
\intertext{Moreover, since $w_1$ is a random $k$-dimensional orthogonal projection of the unit vector $w$, by Corollary~\ref{cor:KDimSphere}, for any $0<\epsilon_2<1$}
\label{eq:UnitBound}
\Prob&\Big(1-\epsilon_2 \leq m\|w_1\|_2^2/k \leq (1-\epsilon_2)^{-1}\Big)\notag\\
& \geq
1-4 e^{-k \cdot c_l(\epsilon_2)/2}.
\end{align}
To conclude, since $\sensing_{\bar{\Lambda}}  \x_{\bar{\Lambda}}$ is Gaussian, its $\ell^2$-norm $\|\sensing_{\bar{\Lambda}}  \x_{\bar{\Lambda}}\|_2^2$ and direction $w$ are mutually independent, hence $\|\sensing_\Lambda^+  w\|_2^2$ and $\|\sensing_{\bar{\Lambda}}  \x_{\bar{\Lambda}}\|_2^2$
are also mutually independent. Therefore,  we can combine the decomposition~\eqref{eq:DeterministicDecomp} with the expected values~\eqref{eq:ExpectFactorA} and~\eqref{eq:ExpectFactorB} to obtain
\begin{align*}
\frac{\Expect \|\Delta_{\textrm{oracle}}(\y,\Lambda) -\x\|_2^2}{\|\x_{\bar{\Lambda}}\|_2^2}
&=
\Expect \|\sensing_\Lambda^+  w\|_2^2
\times
\frac{\Expect \|\sensing_{\bar{\Lambda}}  \x_{\bar{\Lambda}}\|_2^2}{\|\x_{\bar{\Lambda}}\|_2^2} +1\\
&=\frac{k}{m-k-1} +1 
= \frac{1}{1-\frac{k}{m-1}}.\\
\intertext{We conclude that: for any index set $\Lambda$ of size at most $k$, with $k < m-1$, in expectation}
\label{eq:PfOracleError1}
\frac{\Expect \|\Delta_{\textrm{oracle}}(\y,\Lambda) -\x\|_2^2}{\|\x\|_2^2}
&=
\frac{\Expect \|\Delta_{\textrm{oracle}}(\y,\Lambda) -\x\|_2^2}{\|\x_{\bar{\Lambda}}\|_2^2}
\times
\frac{\|\x_{\bar{\Lambda}}\|_2^2}{\|\x\|_2^2}\\
&= \frac{1}{1-\frac{k}{m-1}} \times
\frac{\|\x_{\bar{\Lambda}}\|_2^2}{\|\x\|_2^2}\\
& \geq
\frac{1}{1-\frac{k}{m-1}} \times
\frac{\sigma_k(\x)_2^2}{\|\x\|_2^2}.
\end{align*}

In terms of concentration, combining~\eqref{eq:ConcFactorA},~\eqref{eq:RatioBound}, and~\eqref{eq:UnitBound}, we get that for $0<\epsilon_0,\epsilon_1,\epsilon_2<1$:
\begin{align*}
(1-\epsilon_0)(1-\epsilon_1)(1-\epsilon_2)
\leq & 
\|\sensing^+_\Lambda w\|_2^2
 \frac{\|\sensing_{\bar{\Lambda}}\x_{\bar{\Lambda}}\|_2^2}{\| \x_{\bar{\Lambda}}\|_2^2} \frac{m-k+1}{k}\\
& \leq
[(1-\epsilon_0)(1-\epsilon_1)(1-\epsilon_2)]^{-1}
\end{align*}
except with probability at most (setting $\epsilon_i = \epsilon$, $i=0,1,2$)
\begin{align*}
2 \cdot & e^{-m \cdot c_l(\epsilon_0)/2} + 4 \cdot e^{-k \cdot c_l(\epsilon_2)/2} + 2 \cdot e^{-(m-k+1) \cdot c_l(\epsilon_1)/2}\\
& \leq
8 \cdot e^{-\min(k,m-k+1) \cdot c_l(\epsilon)/2}.
\end{align*}

\subsection{Proof of Theorem~\ref{th:ASCVOracle}}
\label{app:PfThASCVOracle}
Remember that we are considering sequences $k_N,m_N,\sensing_N,\Lambda_N,\x_N$.
Denoting $\slevel_N = k_N/m_N$ and $\mlevel_N = m_N/N$, we observe that the probability~\eqref{eq:OracleKTermConc} can be expressed as $1- 8 e^{-N \cdot c_N(\epsilon)/2}$ where $c_N(\epsilon) = c_l(\epsilon) \cdot \mlevel_N \cdot \min (\slevel_N,1-\slevel_N)$.
For any choice of $\epsilon$, we have
\[
\lim_{N \to \infty} c_N(\epsilon) = c_l(\epsilon) \cdot \mlevel \cdot \min(\slevel,1-\slevel)>0,
\]
hence $\sum_N e^{-N \cdot c_N(\epsilon)/2} < \infty$ and we obtain that for any $\eta>0$
\[
\sum_{N}
\Prob \left(
\left|
\left(\textstyle\frac{\|\Delta_\textrm{oracle}(\y_N,\Lambda_N)-\x_N\|_2^2}{\bestapprox_{k_N}(\x_N)_2^2}-1\right)
\times  \textstyle\frac{m_N-k_N+1}{k_N}
-1
\right|
\geq \eta
\right) < \infty.
\]
This implies \cite[Corollary 4.6.1]{Gray2009} the almost sure convergence
\[
\lim_{N \to \infty} \left(\left(\frac{\|\Delta_\textrm{oracle}(\y_N,\Lambda_N)-\x_N\|_2^2}{\bestapprox_{k_N}(\x_N)_2^2}-1\right)
\times  \frac{m_N-k_N+1}{k_N}\right)
\stackrel{a.s.}{=} 1.
\]
Finally, since $k_N/m_N = \slevel_N \to \slevel$ and $\mlevel_N \to \mlevel$, we also have
\[
\lim_{N\to \infty} \frac{k_N}{m_N-k_N+1} =  \frac \slevel{1-\slevel}
\]
and we conclude that
\begin{align*}
\lim_{N \to \infty} \frac{\|\Delta_{\textrm{oracle}}(\y_N,\Lambda_N) -\x_N\|_2^2}{\|\x_N\|_2^2}
&\stackrel{a.s.}{=} \frac{1}{1-\slevel} \lim_{N \to \infty} \frac{\sigma_{k_N}(\x_N)_2^2}{\|\x_N\|_2^2}\\
& \stackrel{a.s.}{=} \frac{\Gfun{{\pdf}}{2}(\mlevel\slevel)}{1-\slevel}.
\end{align*}

We obtain the result for the least squares decoder by copying the above arguments and starting from~\eqref{eq:ConcentrationIneqLS}.



\subsection{Proof of Lemma~\ref{le:Hbound}}

For the first result we assume that $G(\mlevel^2) \leq (1-\mlevel)^2$. We take $\slevel = \mlevel$ and obtain by definition 
\[
H(\mlevel) \leq \frac{G(\mlevel\slevel)}{1-\slevel} =
\frac{G(\mlevel^2)}{1-\mlevel} \leq (1-\mlevel). 
\]
The second result is a straightforward consequence of the first one.
For the last one, we consider $\mlevel \in (0,\mlevel_0)$. For any
$\slevel\in(0,1)$ we set $\kappa := \mlevel\slevel \in (0,\mlevel_0)$. Since for
any pair $a,b \in (0,1)$ we have $(1-a)(1-b) \leq (1-\sqrt{ab})^2$, we have 
\[
G(\kappa) \geq (1-\sqrt{\kappa})^2 \geq (1-\mlevel)(1-\slevel)
\]
and we conclude that
\[
\forall\ \slevel\in(0,1),\ \frac{G(\mlevel\slevel)}{1-\slevel} \geq 1-\mlevel.
\]

\subsection{Proof of Theorem~\ref{Th:4MomentIntro} and Theorem~\ref{th:4Moment}}

Theorem~\ref{Th:4MomentIntro} and Theorem~\ref{th:4Moment} can be proved from Theorem~\ref{th:ASCVOracle}  and Lemma~\ref{le:Hbound}
along with the following result.

\begin{lemma}\label{le:G4Moment}
Let $\pdf(x)$ be teh PDF of a distribution with finite fourth moment $\Expect X^4 < \infty$. Then there
exists some $\mlevel_0\in(0,1)$ such that the function
$\Gfun{\pdf}{2}(\kappa)$ as defined in Proposition~\ref{prop:RelativeError} satisfies 
\begin{equation}
\Gfun{\pdf}{2}(\kappa) \geq (1-\sqrt{\kappa})^2,\qquad \forall \kappa \in (0,\mlevel_0).\ 
\end{equation}
\end{lemma}

\begin{IEEEproof}[Proof of Lemma~\ref{le:G4Moment}]
Without loss of generality we can assume that $\pdf(x)$ has unit second moment, hence
\[
\Gfun{\pdf}{2}(\kappa) := \frac{\int_0^{\bar{\cdf}^{-1}(1-\kappa)} u^2
 \bar{\pdf}(u) du}{\int_0^\infty u^2 \bar{\pdf}(u) du} =
1-\int_\alpha^\infty u^2 \bar{\pdf}(u) du ,
\]
where we denote $\alpha = \bar{\cdf}^{-1}(1-\kappa)$, which is
equivalent to $\kappa = 1-\bar{\cdf}(\alpha) = \int_\alpha^\infty
\bar{\pdf}(u) du$. 
The inequality $\Gfun{\pdf}{2}(\kappa) \geq (1-\sqrt{\kappa})^2$ is equivalent to
$2\sqrt{\kappa} \geq 1+\kappa-\Gfun{\pdf}{2}(\kappa)$, that is to say 
\begin{equation}
\label{eq:4MomentPf1}
2 \sqrt{\int_\alpha^\infty \bar{\pdf}(u) du}
\geq
\int_\alpha^\infty (u^2+1) \bar{\pdf}(u) du
\end{equation}

By the Cauchy-Schwarz inequality
\[
\int_\alpha^\infty (u^2+1) \bar{\pdf}(u) du
\leq 
\sqrt{\int_\alpha^\infty (u^2+1)^2 \bar{\pdf}(u) du}
\cdot
\sqrt{\int_\alpha^\infty \bar{\pdf}(u) du}.
\]
Since $\Expect X^4 < \infty$, for all small enough $\kappa$
(i.e., large enough $\alpha$), the right hand side is arbitrarily
smaller than $2\sqrt{\int_\alpha^\infty \bar{\pdf}(u) du}$ hence the
inequality $\Gfun{\pdf}{2}(\kappa) \geq (1-\sqrt{\kappa})^2$ holds
true.
\end{IEEEproof}

\begin{IEEEproof}[Proof of Theorem~\ref{Th:4MomentIntro} and Theorem~\ref{th:4Moment}]
Theorem~\ref{Th:4MomentIntro} and Theorem~\ref{th:4Moment} now follow by combining
Lemma~\ref{le:G4Moment} and Lemma~\ref{le:Hbound} to show that 
for a distribution with finite fourth moment there exists a $\mlevel_0
\in (0,1)$ such that $H(\mlevel) \geq 1-\mlevel$ for all $\mlevel \in
(0, \mlevel_0)$. The asymptotic almost sure comparative performance of the estimators then
follows from the concentration bounds in Theorem~\ref{th:ExpectedPerfOracle} and
for the least squares estimator.
\end{IEEEproof}

\subsection{Proof of Proposition~\ref{prop:P0}}

Just as in the proof of Lemma~\ref{le:G4Moment} above, we denote $\alpha = \bar{\cdf}^{-1}(1-\kappa)$, which is equivalent to $\kappa = 1-\bar{\cdf}(\alpha) = \int_\alpha^\infty \bar{\pdf}(u) du$.
We know from Lemma~\ref{le:Hbound} that the identity $\Hfun{\pdf}(\rho) = 1-\rho$ for all $0<\rho<1$ is equivalent to $\Gfun{\pdf}{2}(\kappa) = (1-\sqrt{\kappa})^2$ for all $0<\kappa<1$. 
By the same computations as in the proof of Lemma~\ref{le:G4Moment}, under the unit second moment constraint $\Expect_{\pdf(x)} X^2 = 1$, the latter is equivalent to 
\begin{equation}
\label{eq:4MomentPf1Eq}
2 \sqrt{\int_\alpha^\infty \bar{\pdf}(u) du}
=
\int_\alpha^\infty (u^2+1) \bar{\pdf}(u) du
\end{equation}
Denote $K(\alpha) := \int_\alpha^\infty (u^2+1) \bar{\pdf}(u) du$. The constraint is $K(\alpha) \cdot K(\alpha) = 4 \int_\alpha^\infty \bar{\pdf}(u) du$.
Taking the derivative and negating we must have
\(
2 K(\alpha) \cdot [(\alpha^2+1) \cdot \bar{\pdf}(\alpha)] = 4 \bar{\pdf}(\alpha).
\)
If $\bar{\pdf}(\alpha) \neq 0$ it follows that $K(\alpha) = 2/(\alpha^2+1)$ hence 
\(
(\alpha^2+1) \cdot \bar{\pdf}(\alpha)  =-K'(\alpha) = 4 \alpha / (\alpha^2+1)^2
\)
that is to say
\(
\bar{\pdf}(\alpha) = 4\alpha / (\alpha^2+1)^3
\)
which is satisfied for $\pdf(x) = \pdf_0(x)$.
One can check that
\[
\int_0^\infty \frac{4\alpha}{(\alpha^2+1)^3} d\alpha 
=
\left[ -\frac{1}{(\alpha^2+1)^2}\right]_0^\infty = 1
\]
and, since $\bar{\pdf}(\alpha) \asymp 4 \alpha^{-5}$, $\Expect_{\pdf_0(x)}(X^4) = \infty$.



\subsection{Proof of the statements in Example~\ref{ex:sparseprior}}
Without loss of generality we rescale $\pdf_{\tau,s}(x)$ in the form $\pdf(x) = (1/a) \cdot \pdf_{\tau,s}(x/a)$ so that $p_{\tau,s}$ is a proper PDF with unit variance $\mathbb{E}X^2 = 1$. 
Observing that $\pdf_{\tau,s}(x) \asymp_{x \to \infty} x^{-s}$, we have: $\mathbb{E}X^2 < \infty$ if, and only if $s>3$;  $\mathbb{E}X^4 < \infty$ if, and only if, $s>5$. 
For large $\alpha$, $n=0,2$, $3<s<5$, we obtain 
\[
\int_\alpha^{\infty} x^n \pdf(x) dx 
\asymp \int_\alpha^{\infty} x^{n-s} dx 
\asymp \left[\frac {x^{n+1-s}}{n+1-s}  \right]_\alpha^{\infty} 
\asymp \alpha^{n+1-s}
\]
hence, from the relation between $\kappa$ and $\alpha$, we obtain
\begin{align*}
\frac{1+\kappa-\Gfun{\pdf}{2}(\kappa)}{2\sqrt{\kappa}}
&= \frac{ \int_\alpha^{\infty} (u^2+1) \pdf(u)
 du}{2\sqrt{\int_\alpha^\infty \pdf(u) du}}  
\asymp \frac{
\left( 
\alpha^{3-s} + \alpha^{1-s}
 \right)}{\sqrt {\alpha^{1-s}}} \\
&\asymp \alpha^{\frac{5-s}{2}}
\end{align*}
For $3 < s < 5$ we get
\[
\lim_{\kappa \to 0} \frac{1+\kappa-\Gfun{\pdf}{2}(\kappa)}{2\sqrt{\kappa}} = \infty
\]
hence there exists $\delta_0>0$ such that for $\kappa < \sqrt\delta_0$
\[
\Gfun{\pdf}{2}(\kappa) < 1+\kappa-2\sqrt\kappa =  (1-\sqrt{\kappa})^2.
\]
We conclude using Lemma~\ref{le:Hbound}.

\subsection{The Laplace distribution}
\label{app:Laplace}
First we compute $\bar{\pdf}_1(x) = \exp(-x)$ for $x\geq 0$,
$\bar{\cdf}_1(z) = 1-e^{-z}$, $z \geq 0$ hence
$\bar{\cdf}_1^{-1}(1-\kappa) = -\ln \kappa$.  
For all integers $q \geq 1$ and $x >0$, we obtain by integration by
parts the recurrence relation 
\[
\int_0^x u^q e^{-u} du = q \int_0^x u^{q-1} e^{-u} du - x^q e^{-x},
\forall q\geq 1. 
\]
$
\int_0^x e^{-u} du = 1-e^{-x},
$
hence for $q = 1$ we obtain
$
\int_0^x u e^{-u} du = 1-e^{-x} - x e^{-x} = 1- (1+x) e^{-x},
$
and for $q=2$ it is easy to compute
\[
\int_0^x u^2 e^{-u} du 
= 2-(2+2x+x^2) e^{-x}
\]
\eqref{eq:BestKTermLaplaceAsymp} and
\eqref{eq:BestKTermLaplaceAsympL1} follow from substituting these
expressions into:
\[
\Gfun{\pdf_1}{q}(\kappa) = 
\frac{\int_0^{-\ln \kappa} u^q \bar{\pdf}_1(u) du}{\int_0^\infty u^q
  \bar{\pdf}_1(u) du}.
\]

\bibliographystyle{plain}
\bibliography{main}

\begin{IEEEbiographynophoto}{R. Gribonval}
  graduated from {\'E}cole Normale Sup{\'e}rieure, Paris, France. He
  received the Ph. D. degree in applied mathematics from the University of
  Paris-IX Dauphine, Paris, France, in 1999, and his Habilitation {\`a} Diriger des Recherches in applied mathematics from the University of Rennes I, Rennes, France, in 2007. 
    He was a  visiting scholar at the Industrial Mathematics Institute, University of South Carolina during 1999-2001. He is now a Directeur de Recherche with Inria in Rennes, France, which he joined in 2000. His research interests range from mathematical signal processing to machine learning and their applications, with an emphasis on multichannel audio and compressed sensing. He founded the series of international workshops SPARS on Signal Processing with Adaptive/Sparse Representations. In 2011, he was awarded the Blaise Pascal Award in Applied Mathematics and Scientific Engineering from the French National Academy of Sciences, and a Starting Grant from the European Research Council. He is a Senior Member of the IEEE since 2008 and a member of the IEEE Technical Committee on Signal Processing Theory and Methods since 2011.
\end{IEEEbiographynophoto} 

\begin{IEEEbiographynophoto}{Volkan Cevher}
received his BSc degree (valedictorian) in Electrical Engineering from Bilkent University in 1999, and his PhD degree in Electrical and Computer Engineering from Georgia Institute of Technology in 2005. He held research scientist positions at University of Maryland, College Park during 2006-2007 and at Rice University during 2008-2009. Currently, he is an assistant professor at Swiss Federal Institute of Technology Lausanne and a Faculty Fellow at Rice University. His research interests include signal processing theory and methods, machine learning, and information theory. Dr. Cevher received a best paper award at SPARS in 2009 and an ERC StG in 2011. 
\end{IEEEbiographynophoto}

\begin{IEEEbiographynophoto}{Mike E. Davies} (M'00-SM'11) received the B.A. (Hons.) degree in engineering from Cambridge University, Cambridge, U.K., in 1989 and the Ph.D. degree in nonlinear dynamics and signal processing from University College London (UCL), London, U.K., in 1993. Mike Davies was awarded a Royal Society University Research Fellowship in 1993 which he held first at UCL and then the University of Cambridge. He acted as an Associate Editor for IEEE Transactions in Speech, Language and Audio Processing, 2003-2007.

Since 2006 Dr. Davies has been with the University of Edinburgh where he holds the Jeffrey Collins SFC funded chair in Signal and Image Processing. He is the Director of the Joint Research Institute in Signal and Image Processing, a joint venture between the University of Edinburgh and Heriot-Watt university as part of the Edinburgh Research Partnership. His current research focus is on sparse approximation, computational harmonic analysis, compressed sensing and their applications within signal processing. His other research interests include: non-Gaussian signal processing, high-dimensional statistics and information theory.

\end{IEEEbiographynophoto}
\vfill
\end{document}